\documentclass{amsart}

\usepackage{mathtools, amssymb, hyperref, amsmath,amssymb,amsbsy,amsfonts, latexsym,amsopn,amstext, amsxtra,euscript,amscd, amsthm}
\usepackage{amssymb, amsthm, amsmath,  amssymb, amsbsy,   amsfonts,  latexsym,  amsopn,   amstext, amsxtra,  euscript,   amscd}
\usepackage{amsaddr}
\usepackage{tikz}
\usepackage[capitalise]{cleveref}
\usepackage{booktabs, algorithm, algorithmic, graphicx}
\usepackage[msc-links, backrefs]{amsrefs}

\usepackage{enumitem}

\hypersetup{
  colorlinks   = true, 
  urlcolor     = blue, 
  linkcolor    = blue, 
  citecolor   = red 
}

\newtheorem{thm}{Theorem} 
 
\newtheorem{lem}{Lemma}

\newtheorem{exa}{Example}

\newtheorem{rem}{Remark}
\newtheorem{cor}{Corollary}

 \crefname{thm}{Thm.}{}
\crefname{prop}{Prop.}{}
\crefname{lem}{Lem.}{}
\crefname{cor}{Cor.}{}
\crefname{figure}{Fig.}{}

\newcommand\Aut{Aut}
\newcommand\h{h}

\newcommand\Z{\mathbb Z}
\newcommand\Q{\mathbb Q}
\newcommand\R{\mathbb R}
\newcommand\C{\mathbb C}
\newcommand\bQ{\overline{\mathbb Q}}
 
\newcommand\M{{\mathcal M}}
\newcommand\X{\mathcal X}     
\newcommand\B{\mathcal B}

\newcommand\Jac{\mbox{Jac }}

\newcommand\GL{GL}
\newcommand\chara{char}

\newcommand{\wgcd}{\mathrm{wgcd}}
\newcommand{\awgcd}{\overline{\mathrm{wgcd}}}

\newcommand{\wh}{\mathfrak{H}_k}
\newcommand{\awh}{\mathfrak{H}}

\newcommand{\alwh}{\mathfrak{s}}

\newcommand\Img{Img }

\newcommand\Hom{Hom }
\newcommand\End{End }

\newcommand{\wP}{\mathbb{P}}                          	
\newcommand\x{\mathbf x}
\newcommand\y{\mathbf y}
\newcommand\A{\mathbb A}

\newcommand\w{\mathfrak w}

\newcommand\CC{\mathcal C}
\newcommand\K{\mathcal K}

\newcommand\p{\mathfrak p}

\begin{document}

\title[Machine learning for moduli space of genus two curves]{Machine learning for moduli space of genus two curves and an application to isogeny based cryptography}
 
\author{Elira Shaska}
\address{Department of Computer Science \\ 
College of Computer Science and Engineering, \\
Oakland  University,  Rochester, MI, 48309. }
\email{elirashaska@oakland.edu}

\author{Tanush Shaska}
\address{Department of Mathematics and Statistics, \\
College of Arts and Sciences \\
 Oakland University, Rochester, MI, 48309}
\email{tanush@umich.edu}

\begin{abstract}
We use machine learning to study the moduli space of genus two curves,  specifically  focusing on 
detecting whether a genus two curve has $(n, n)$-split Jacobian.  
Based on such techniques, we observe that  there are very few rational moduli points with small weighted moduli height and $(n, n)$-split Jacobian for $n=2, 3, 5$. We computational prove that 
there are only 34 genus two curves (resp. 44 curves)    with (2,2)-split Jacobians (resp. (3,3)-split Jacobians)   and weighted moduli height $\leq 3$.  

We discuss different machine learning models for such applications and demonstrate the ability to detect splitting with high accuracy using only the Igusa invariants of the curve. This shows that artificial neural networks and machine learning techniques can be highly reliable for arithmetic questions in the moduli space of genus two curves and may have potential applications in isogeny-based cryptography.
\end{abstract}

\keywords{genus two curves, isogenies, weighted projective space}
 
\subjclass[MSC 2020]{ Primary:   68T07;    Secondary:  68T20, 68Q32}

\maketitle
 
\section{Introduction}
This paper is a first  attempt to use machine learning to study the moduli space $\M_g$  of algebraic curves of genus $g \geq 2$.  We focus on the simplest case when the genus is $g=2$ and curves are defined over   $\Q$, since we are interested in arithmetic properties of the moduli space $\M_2$,  rational points of  $\M_2$,  and in particular on fine moduli points (i.e, those rational points for which exists  a curve $C$ defined over the rationals). 

We create a database of   isomorphism classes of genus two curves defined over $\Q$  which  are uniquely determined  by points in the weighted projective space $\wP_{(2,4,6,10)}$  given by   classical  invariants of binary sextics  $J_2, J_4, J_6, J_{10}$.    
To determine uniquely a point in $\wP_{(2,4,6,10)}$  we normalize    each point $\p$ by  the absolute weighted greatest common divisor 
$\lambda:={\awgcd(\p)}$,   by multiplying    $\frac 1 \lambda \star \p$. This assures that these moduli points have coordinates as small as possible.
%
The data is stored in a Python dictionary where the key is the triple of absolute invariants $(t_1, t_2, t_3)$ as defined in \cref{t-invariants}.  Since the these triples uniquely determine the moduli points, this assures that there is no redundancy in our data.        In all our methods the keys $(t_1, t_2, t_3)$ are not used in any training, but only for graphical purposes.  The training is done only with  points $\p \in \wP_{(2,4,6,10)}$. 

Most of this paper is on creating the training dataset, which is  based mostly on methods described in previous work of the authors; 
\cite{MS-2, obus-sh, cms-1, cms-2, cms-3, sha-2016}.  
 Data points are ordered by their weighted moduli height as defined in  \cite{s-sh}.  For every point $\p \in \wP_{(2,4,6,10)}$ we compute the absolute weighted moduli height $\awh (\p)$  of $\p$, its automorphism group $\Aut (\mathfrak p)$, the conic $Q:=\p/ \Aut (\mathfrak p)$ and determine  
 whether the point is a fine moduli point (i.e, it is defined over its field of moduli) or a coarse moduli point depending whether the conic $Q$ has or not a rational point.  

Our first problem is to determine the distribution of fine moduli points in the moduli space and to see what happens to the number of fine points as the weighted moduli height grows to infinity.    

The second problem that is considered is to use machine learning to determine  if a point $\p \in \wP_{(2,4,6,10)}$ has split Jacobian. 
Genus two curves with $(n, n)$-split Jacobians have been the focus of study for a long time starting with Jacobi, Hermite, Goursat, 
Fricke, Brioschi, and in the last few decades by Frey, Kani, Shaska, Kumar, et al.   For $n$ odd, the locus of genus 2 curves with $(n, n)$-split Jacobian, denoted by Shaska in \cite{sha-2001, thesis} by ${\mathcal L}_n$, is a 2-dimensional irreducible locus in the moduli space of genus 2 curves $\M_2$.   
When $n=2$  the covering is Galois and that implies that the genus two curve has an extra involution. Such curves were studied from Jacobi and many authors since; see \cite{deg2} for a modern treatment.  

Let us now consider only the cases when $n$ is odd. The equation of $\mathcal L_n$ were first computed by Shaska  in his thesis (2001)  in  \cite{thesis} and  papers published in   \cite{deg-3, deg5}.  A  general method was given for all $n$ in \cite{sha-2001} and \cite{thesis}.  In 2015,  Kumar   confirmed such computational via a different approach; see  \cite{kumar}. 

With the recent developments on isogeny based cryptography such curves have received new attention.  There are many papers on trying to decide if  a genus two curve has $(n, n)$-split Jacobian for small $n$.   In   \cite{costello-nn-split} was shown that detecting if a Jacobian is $(n,n)$-split can speed computation of isogenies considerably.  How do you detect of a genus two curve is in $\mathcal L_n$?  For small $n$ we can just plug the invariants of the curve in the equation of $\mathcal L_n$, even though  many authors seem unaware of how to obtain such equations.  For larger $n$ computing such equations is impractical which suggests a machine learning method could be much more feasible. 
  Since ${\mathcal L}_n$ are rational surfaces, there must be a way to  generate rational points in ${\mathcal L}_n$.     Hence, even without knowing explicitly the equation of ${\mathcal L}_n$, we can create a training dataset and use such data to train some machine learning model. 

After experimenting with many methods of unsupervised and supervised learning we were able to detect the splitting using K-neighbours Classifier with an accuracy of 99.94\% when $n=2, 3, 5$.  There is no reason to believe that this accuracy will go down for higher $n$.  
Our experiments were run for randomly generated dataset of 50 000 points for each locus $\mathcal L_n$, $n=2, 3, 5$.  It seems as with more computational resources such results can be easily extended to larger $n$.   It remains to be seen how such experiments can be performed over a field of positive characteristic.

On the first problem: distribution of fine moduli points we were less successful. The reason was that there are very few fine points for moduli points of weighted height $\leq 3$.  We can generate randomly such points, using for example the locus $\mathcal L_2$, but they would be rather useless if we are trying to study the distribution of fine points in $\mathcal M_2$.  Moreover, we are more interested in fine moduli points which are in $\mathcal M_2 \setminus \mathcal L_2$.  In order to do  meaningful experiments it seems that we need to have a larger dataset (say up to weighted height $\wh \leq 5$) which requires computational resources which are not available to us.  Due to these experiments, we were able to discover some interesting results. For example we showed that 
there are no rational points $\mathfrak p\in {\mathcal L}_2$ with weighted  height $\wh (\p) < 3/2$, or that there are exactly 34 such points with $\wh (\p) < 3$.   For $\mathcal L_3$ we showed that there are no rational points with weighted height $\wh <2$ and there are 46 such points with $2 \leq \wh (\mathfrak p) <3$. Such results would have been difficult to notice without machine learning techniques, which highlights exactly the way  Machine Learning is the most likely to be used in mathematics: detecting a pattern and then trying to prove it by brute force or other classical methods.

\section{Preliminaries}\label{sec-2}
Since this is intended for non experts  in machine learning we  briefly give a basic setup. 
We also provide basic definitions of weighted projective spaces and weighted heights, since all our computations are done in the weighted space.

\subsection{Neural networks}
Let $k$ be a field. 
A \textbf{neuron} is a function $f : k^n \to k$ such that  for every $\x = (x_1, \ldots , x_n) $     we have 
\[
f (\x) = \sum_{i=1}^n  w_i x_i  + b,
\]
where $b\in k$ is a constant   called \textbf{bias}.    We can generalize neurons to  tuples of neurons via   $F:  \;  k^n   \to k^n$
\[
\x   \to  \left( f_1   (\x), \ldots , f_n (\x)  \right)
\]
for any given set of positive integers $\w_0, \ldots , \w_n$ called  weights.  
Then $F$  is a   function    written as 
$
F (\x )=  W \cdot \x + \beta,
$
for some $\beta \in k^{n}$ and  $W$  an $n \times n$ matrix with integer entries. 
A  (non-linear) function  $g: k^n  \to k^n $   is called an \textbf{activation function} while   a \textbf{network layer} is a function   $L: k^n \to k^n$, such that. 
\[
\x   \to g \left( W\cdot \x + \beta  \right)
\]
 for some  activation function $g$.   A \textbf{neural network} is the composition of many layers.  The $i$-th layer 
\[
\begin{split}
\cdots \longrightarrow k^n    &  \stackrel{L_i}{\longrightarrow} k^n  \longrightarrow \cdots  \\
\x & \longrightarrow  L_i  (\x) =  g_i  \left( W^i  \x + \beta^i  \right),
\end{split}
\]
where $g_i$, $W^i$, and $\beta^i$ are the activation, matrix, and bias  for   this layer. 
  
\subsection{Weighted greatest common divisors}
Let $\x = (x_0, \dots x_n ) \in \Z^{n+1}$ be a tuple of integers, not all equal to zero.  
A set of \textbf{weights} is  the ordered tuple   $ \w=(q_0, \dots, q_n)$, where $q_0$, \dots, $q_n$ are positive integers.
A \textit{weighted integer tuple} is a tuple $\x = (x_0, \dots, x_n ) \in \Z^{n+1}$ such that to each coordinate $x_i$ is assigned the weight $q_i$.  We multiply weighted tuples by scalars $\lambda \in \Q$ via 
\[ \lambda \star (x_0, \dots , x_n) = \left( \lambda^{q_0} x_0, \dots , \lambda^{q_n} x_n   \right) \]
For an ordered   tuple of integers  $\x=(x_0, \dots, x_n) \in \Z^{n+1}$, whose coordinates are not all zero, the \textbf{weighted greatest common divisor with respect to} $\w=(q_0, \ldots , q_n)$ is the largest integer $d$ such that 
\[ d^{q_i} \, \mid \, x_i, \; \; \text{for all }    i=0, \dots, n.\]
We will call a point $\p \in {\mathbb P}_{\w}^n (\Q)$ \textbf{normalized} if $\wgcd (\p) =1$. 
The \textbf{absolute weighted greatest common divisor} of  an integer tuple $\x=(x_0, \dots , x_n)$ with respect to the set of weights $\w=(q_0, \dots , q_n)$ is the largest  real number $d$ such that %
\[ d^{q_i} \in \Z  \quad {\text and} \quad   d^{q_i}\, \mid \, x_i, \; \; \text{for all }    i=0, \dots n.\]
It is denoted by $\awgcd (\x)$.    We will call a point $\p \in {\mathbb P}_{\w}^n (\Q)$ \textbf{absolutely normalized} if $\awgcd (\p) =1$. 

\subsection{Weighted projective spaces and moduli space of binary forms}
%
Consider the action of   $k^\ast = k \setminus \{0\}$ on $\A_k^{n+1} \setminus \{(0, \cdots, 0)\}$ given by
\begin{equation}\label{equivalence}
\lambda \star (x_0, \dots , x_n) = \left( \lambda^{q_0} x_0, \dots , \lambda^{q_n} x_n   \right), \      \text{for }\ \lambda \in k^\ast.
\end{equation}
Define the \textbf{weighted projective space}, denoted by      $\wP_\w^n (k)$,   to be  the quotient space  $\A_k^{n+1}/k^\ast$ of this action, which is  a geometric quotient since $k^\ast$ is a reductive group.    An element    $\x \in \wP_\w^n (k)$ is denoted by $\x = [ x_0 :   \dots : x_n]$  and its $i$-th coordinate  by $x_i (\x)$.  For more details on weighted projective spaces we refer to \cites{Beltrametti1986, s-sh}.

\subsubsection{Heights on weighted projective spaces}
 Heights on weighted projective spaces can be defined for any number field $k$ and its integer ring $\mathcal O_k$.  We will state the definitions over a number field $k$ even though for our computations we will only consider the case $k=\Q$.

Let $\w=(q_0, \dots , q_n)$ be a set of weights and ${{\mathbb P}_{\w}^n} $ the weighted  projective space  over  a number field $k$.   Let  $\p=[x_0: \dots : x_n]$ be a point in 
${{\mathbb P}_{\w}^n} (k)$.  Without loss of generality we can assume that $\p =[x_0 : \dots : x_n] $   has coordinates in $(\mathcal O_k$. Let $M_k$ be the set of places of $k$.

We follow  definitions in \cite{b-g-sh} and define   \textbf{weighted multiplicative height} of $\p$   by
\begin{equation}\label{def:height}
\wh( \p ) := \prod_{v \in M_k} \max   \left\{   \frac{}{}   |x_0|_v^{\frac {n_v} {q_0}} , \dots, |x_n|_v^{\frac {n_v} {q_n}} \right\}.
\end{equation}
The \textbf{absolute   weighted height} of $\p \in {{\mathbb P}_{\w}^n} (k)$  is the function   $\awh: {{\mathbb P}_{\w, \bQ}^n}   \to [1, \infty)$, 
\[
\awh  (\p)=\wh(\p)^{1/[k:\Q]},
\]
%
for any $k$ which contains $\Q (\awgcd (\p))$; see \cite{b-g-sh} and \cite{s-sh} for details. 

Let ${\mathbb P}_{\w, k}$ be a well-formed weighted projective space and   $\x= [ x_0 : \dots : x_n]  \in {\mathbb P}_{\w, k} (k)$.  Assume $\x$ normalized (i.e. $\wgcd_k (\x)=1$).    Clearly $\wgcd (\x) |  \gcd (x_0, \ldots , x_n)$  and therefore $\wgcd (\x) \leq \gcd (x_0, \ldots , x_n)$.
Let      $\x$ be absolutely normalized. Then 
$\gcd (x_0, \ldots , x_n)=1$.  
If  $\x=[x_0 : \cdots : x_n]$ is a normalized point then by definition of the height
\[
\wh (\x) = \max_{i=0}^n \{   | x_i |^{\frac 1 {q_i}} \}
\]
Assume now that $\x= [ \lambda^{q_0} x_0 : \cdots : \lambda^{q_n} x_n ]$ such that
$
\lambda = \wgcd \left(  \lambda^{q_0} x_0 : \cdots : \lambda^{q_n} x_n  \right).
$
Denote by $s$ the index where $\min_j  \{ |  \lambda^{q_i} \, x_j|^{\frac 1 {q_j}} \} = \lambda  \min_j  \{ |   x_j|^{\frac 1 {q_j}} \}  $ is obtained.
Then
\begin{equation}\label{eq-1}
\frac 1 {\lambda \cdot  x_s^{1/q_s}} \star \x = 
\left[
\frac {x_0} {x_s^{q_0/q_s}}: \cdots  : 1 : \cdots :  \frac {x_n} {x_s^{q_n/q_s}}
\right]
= : \y 
\end{equation}
where 1 is in the $s$ position. Simplify all coordinates in \cref{eq-1}.  Multiplying $\y$ by $x_s^{ \frac 1 {q_s} } $ we have 
$
x_s^{ \frac 1 {q_s} }   \star \y = [ x_0: \cdots : x_n],
$
which is now a normalized point.    Hence
\[
\wh (\x) = \frac 
{ \max_{i=0}^n \{    |\lambda^{q_i} x_i |^{ \frac 1 {q_i}}  \} }  
{\min_{i=0}^n \{    |\lambda^{q_i} x_i |^{ \frac 1 {q_i}} \}   }
\]
Notice that    $\wh (\x)$ is given by  $ \wh (\x) =      \frac   {\max_i  |x_i|^{\frac 1 {q_i}}    }     {\min_j  | x_j|^{\frac 1 {q_j}} }$.
When   $\x=[ x_0 : \dots : x_n] \in {\mathbb P}_{\w} (\Q)$ is an absolutely    normalized point, then    $\gcd (x_0, \ldots , x_n) \le \wh (\x)$. 

In \cref{sec-data} we create a database of isomorphism classes of genus two curves as a set of points in a weighted projective space. For each point in that database we compute the weighted height and the absolute weighted height. 

\section{Abelian surfaces and their  isogenies}
Here we provide the minimum  setup of isogenies of Abelian varieties and their use in cryptography.  The interested reader should check  \cite{frey-sh} among many other places for further details.

Let ${\mathcal A}$, $\B$ be abelian varieties over a    field $k$.  We denote the $\Z$-module of homomorphisms  ${\mathcal A} \mapsto  \B$  by $\Hom( {\mathcal A}, \B)$  and the ring of endomorphisms ${\mathcal A} \mapsto {\mathcal A}$ by $\End {\mathcal A}$. It is more convenient  to   work with the $\Q$-vector spaces 
\[
\Hom^0 ({\mathcal A}, \B):= \Hom({\mathcal A}, \B) \otimes_\Z \Q,
\]
 and $\End^0 {\mathcal A}:= \End {\mathcal A}\otimes_\Z \Q$.    
 A homomorphism $f : {\mathcal A} \to  \B$ is called an \textbf{isogeny} if $\Img (f) = \B$  and $\ker (f)$ is a finite group scheme. If an isogeny ${\mathcal A} \to \B$ exists we say that ${\mathcal A}$ and $\B$ are isogenous.  
The degree of an isogeny $f : {\mathcal A} \to \B$ is the degree of the function field extension
\[ \deg f  := [k ({\mathcal A}) : f^\star k(\B)].\]
It is equal to the order of the group scheme $\ker (f)$.     
The group of $\bar{k}$-rational points has order 
\[
\#(\ker f)(\bar{k}) = [k(A) : f^\star k(B)]^{sep},
\]
 where $[k(A) : f^\star k(B)]^{sep}$ is the degree of the maximally separable extension in $k({\mathcal A})/ f^\star k(\B)$.   $f$ is a \textbf{separable isogeny} if and only if 
$ \# \ker f(\bar{k}) = \deg f.$
The   basic principle on any isogeny based cryptosystem is based on the fact that 
for any Abelian variety ${\mathcal A}/k$ there is a one to one correspondence between the finite subgroup schemes $\K \leq {\mathcal A}$ and   isogenies $f : {\mathcal A} \to \B$, where $\B$ is determined up to isomorphism.   Moreover, $\K = \ker f$ and $\B = {\mathcal A}/\K$.

If  ${\mathcal A}$ and $\B$ are isogenous then $\End^0 ({\mathcal A}) \cong \End^0 (\B)$.     The following is often called the fundamental theorem of Abelian varieties:

\begin{thm}[Poincare-Weil]
Let ${\mathcal A}$ be an Abelian variety.  Then ${\mathcal A}$ is isogenous to 
\[   {\mathcal A}_1^{n_1} \times {\mathcal A}_2^{n_2} \times \dots \times {\mathcal A}_r^{n_r}, \]
where (up to permutation of the factors) ${\mathcal A}_i$ , for $i=1, \dots , r$ are simple, non-isogenous, Abelian varieties. 
Moreover,  up to permutations,  the factors  ${\mathcal A}_i^{n_i}$ are uniquely determined  up to isogenies. 
\end{thm}

If ${\mathcal A}$ is a  absolutely simple Abelian variety then every endomorphism  not equal $0$ is an isogeny. 
 
\subsubsection{Computing isogenies between Abelian varieties} 

Fix a field $k$ and let ${\mathcal A}$ be an Abelian variety over $k$.  Let $H$ denote a finite subgroup of ${\mathcal A}$.  From the computational point of view we have the following problems:

\begin{enumerate}
\item Compute all Abelian varieties $\B$ over $k$ such that there exists an isogeny ${\mathcal A} \to \B$ whose kernel is isomorphic to $H$. 

\item Given ${\mathcal A}$ and $H$, determine   $\B:={\mathcal A}/H$ and the isogeny ${\mathcal A} \to \B$.

\item Given two Abelian varieties ${\mathcal A}$ and $\B$, determine if they are isogenous and compute a rational expression for an isogeny ${\mathcal A} \to \B$. 
\end{enumerate}
 
\subsection{Torsion points and dual isogenies}
The most classical example of an isogeny is the  scalar multiplication by $n$ map  
$
[n] : \;  {\mathcal A} \to {\mathcal A}.
$
The kernel of $[n]$ is a group scheme of order $n^{2\dim {\mathcal A}}$.  We denote  by ${\mathcal A} [n]$  the group $\ker  [n] (\bar{k}) $.   The elements in ${\mathcal A}[n]$   are called $n$-\textbf{torsion points} of ${\mathcal A}$.

Let $f : {\mathcal A} \to \B$ be a degree $n$ isogeny.  Then there exists an isogeny $\hat f : \B \to {\mathcal A}$ such that
\[ f \circ \hat f = \hat f \circ f = [n]. \]
The isogeny $\hat f$ is called the \textbf{dual } of $f$.
 
\begin{thm}\label{thm-1} 
Let ${\mathcal A}/k$ be an Abelian variety, $p = \chara k$, and $\dim {\mathcal A}= g$. 

\begin{itemize}
\item[i)] If $p \nmid \, n$, then $[n]$ is separable, $\# {\mathcal A}[n]= n^{2g}$ and ${\mathcal A}[n]\cong (\Z/n\Z)^{2g}$.

\item[ii)] If $p \mid n$, then $[n]$ is inseparable.  Moreover, there is an integer $0 \leq i \leq g$ such that 
\[ {\mathcal A} [p^m] \cong (\Z/p^m\Z)^i, \; \text{for all } \; m \geq 1. \]
\end{itemize}
\end{thm}

If $i=g$ then ${\mathcal A}$ is called \textbf{ordinary}.  If ${\mathcal A}[p^s](\bar k)= \Z/p^{ts}\Z$    
then the abelian variety has \textbf{$p$-rank} $t$. If $\dim {\mathcal A}=1$ (elliptic curve) then it is called \textbf{supersingular} if it has $p$-rank 0.
An abelian variety ${\mathcal A}$ is called \textbf{supersingular} if it is isogenous to a product of supersingular elliptic curves.  

\subsection{Supersingular Isogeny Diffie-Helman (SIDH)} 
Let ${\mathcal A}= \Jac (\CC)$, where $\CC$ is a genus     $g\geq 2$  curve  defined over a field $k$ with $\chara k = p$, for $p\neq 2$.  
 We pick $p=m\cdot n -1$ where $(m, n)=1$.  
 Then ${\mathcal A} [m]$ and ${\mathcal A} [n]$ have dimension $2g$.  Pick bases 
\[
 \{P_1 , \ldots , P_{2g} \} \;  \text{ for }   {\mathcal A} [m]  \quad \text{  and } \quad  \{Q_1 , \ldots , Q_{2g} \} \;  \text{ for } \;  {\mathcal A} [n]  
 \]
Alice randomly generates 
\[
   R_1     =\sum_{i =1}^{2g} [x_{1, i}] P_i,  \; 
   R_2  =\sum_{i=1}^{2g} [x_{2, i}] P_i, \;    \ldots ,  \; 
    R_{2g-1}  =\sum_{i=1}^{2g} [x_{2g-1, i}] P_i, 
\]
and a maximal subgroup $\K_A = \langle R_1, \ldots , R_{2g-1} \rangle$.   Then she computes the isogeny  $\phi_A : {\mathcal A} \to {\mathcal A}/ \K_A =: {\mathcal A}_A$.   Alice generates  a public key  $( {\mathcal A}_A, \phi_A (Q_1), \ldots , \phi_A (Q_{2g} )$.     Bob randomly generates 
\[
    S_1=\sum_{i =1}^{2g} [x_{1, i}] Q_i,   \; 
  S_2=\sum_{i=1}^{2g} [x_{2, i}] Q_i, \; 
   \ldots    \; 
  S_{2g-1}=\sum_{i=1}^{2g} [x_{2g-1, i}] Q_i    
\]
and a maximal subgroup $\K_B = \langle S_1, \ldots , S_{2g-1} \rangle$.   Then he computes the isogeny     \[\phi_B : {\mathcal A} \to {\mathcal A} / \K_B =: {\mathcal A}_B.\]
Bob generates  a public key  $( {\mathcal A}_B, \phi_B (P_1), \ldots , \phi_B (P_{2g} )$.  

\subsubsection{Key exchange:}   Using   $\phi_B (P_1), \ldots , \phi_B (P_{2g})$, 
Alice generates 
\[
    R_1^\prime =\sum_{i =1}^{2g} [x_{1, i}] \phi_B (P_i),  \;  
  R_2^\prime =\sum_{i=1}^{2g} [x_{2, i}] \phi_B (P_i), \;     \ldots  \; 
  R_{2g-1}^\prime =\sum_{i=1}^{2g} [x_{2g-1, i}] \phi_B (P_i)
\]
and creates $\K_A^\prime = \langle R_1^\prime, \ldots, R_{2g-1}^\prime \rangle$, a maximal subgroup of ${\mathcal A}_B [n]$ and the isogeny 
\[
\phi_A^\prime : {\mathcal A}_B \to {\mathcal A}_B/ \K_A^\prime = : {\mathcal A}_{BA}.
\]
Using   $\phi_A (Q_1), \ldots , \phi_A (Q_{2g}$,     Bob generates 
\[
  S_1^\prime =\sum_{i =1}^{2g} [x_{1, i}] \phi_A (Q_i), \;
   S_2^\prime =\sum_{i=1}^{2g} [x_{2, i}] \phi_A (Q_i), \;   \ldots ,  \; 
  S_{2g-1}^\prime =\sum_{i=1}^{2g} [x_{2g-1, i}] \phi_A (Q_i), 
\]
and creates $\K_B^\prime = \langle S_1^\prime, \ldots, S_{2g-1}^\prime \rangle$, a maximal subgroup of ${\mathcal A}_A [m]$ and the isogeny 
\[
\phi_B^\prime : {\mathcal A}_A \to {\mathcal A}_A/ \K_B^\prime = : {\mathcal A}_{AB}.
\]
\begin{lem}\label{isog}
We have that $\phi_B^\prime  \circ   \phi_A   =  \phi_A^\prime  \circ \phi_B$. 
Moreover,    ${\mathcal A}_{AB} = {\mathcal A}_{BA}$.
\end{lem}

Let  $\CC$ be a genus 2 curve and $\psi : \CC \to E_1$ be a maximal degree $n$ covering which does not factor through an isogeny. Then, there is another elliptic curve $E_2 := \Jac \CC / E_1$ such that $\Jac \CC$ is isogenous via a degree $n^2$ isogeny to the product $E_1 \times E_2$.  We say that $\Jac \CC$ is \textbf{$(n, n)$-decomposable} or \textbf{$(n, n)$-split}.  
The locus of curves with $(n, n)$-decomposable Jacobians is a 2-dimensional irreducible locus in $\M_2$. Such loci for small $n=2, 3, 5$ are computed in 
\cite{sha-2001, thesis, deg2, deg-3,  deg5}.    

$\Jac (\CC)$ is a geometrically simple Abelian variety if and only if it is not $(n, n)$-split for some $n>1$.  In other words, if 
$\Jac (\CC)$ is split over $k$, then there exists an integer $n\geq 2$ such that $\Jac (\CC)$ is $(n, n)$-split. 

\section{Moduli space $\M_2$ of genus 2 curves}\label{sec-4}
We refer to \cite{rat} for notation and terminology. We assume $k$ is a field of $\mbox{char} k \neq 2$, 
so a genus two curve $C$  is given by an  equation 
\[
C : \qquad z^2 y^4 = a_6 x^6 + a_5 x^5  y+ \cdots a_1 x y^5 + a_0 y^6
\]
such that the discriminant of the sextic on the right is nonzero.  Hence, the isomorphism class of  $C$ is determined by the invariants of $f(x, y)$, which are commonly denoted by $J_2(f), J_4(f), J_6(f), J_{10} (f)$ and are homogenous polynomials of degree 2, 4, 6, and 10 respectively in the coefficients of $f( x, y)$.
Moreover, the invariant $J_{10} (f)$   is the discriminant of the sextic and therefore $J_{10} (f) \neq 0$.   Hence, the  moduli space of genus 2 curves defined over $k$ is isomorphic to $\wP_{(2,4,6,10)} \setminus \{J_{10} \neq 0\}$. 
There is a   morphism   $\wP_{(2,4,6,10)}^3          \to \wP_{(1,2,3,5)}^3,$   given by 
\begin{equation}\label{iso-wp}
[x_0: x_1:x_2:x_3]     \to [y_0: y_1 :y_2:y_3]=\left[  x_0^2 : x_1^2 : x_2^2: x_3^2        \right]
\end{equation}
Since we want to design a model where the incoming features will be a genus two curve, then equivalently this means a point 
$[x_0: x_1:x_2:x_3]  \in \wP_{(1,2,3,5)}$.  Equivalently the input could be the equation of the curve, but this poses no issue for $g=2$ since we can  easily compute invariants.   However, even though the space $\wP_{(1,2,3,5)}$ is nicer to deal with, from the computational point of view we rather store points in the form 
$(J_2, J_4, J_6, J_{10})$ than $(J_2^2, J_4^2, J_6^2, J_{10}^2)$.  Hence, our database will consists of points of $\wP_{(2,4,6,10)}$.

There is also an issue to address when it comes to finding the "smallest" representatives for the equivalence class $[x_0: x_1:x_2:x_3]$.  Theoretically this is handled in \cite{b-g-sh}, but that would require computing weighted greatest common divisors and that could be very costly for large coordinates $x_0, \ldots , x_3$. 
Since $q=1\cdot 2\cdot 3 \cdot 5= 30$,   the  Veronese map is 
\begin{equation}
\left[  J_2: J_4: J_6: J_{10} \right] \longrightarrow \left[  J_2^{30}: J_4^{15}: J_6^{10}: J_{10}^6 \right]
= \left[
\frac {J_2^{30}} {J_{10}^6 } :  \frac {J_4^{15}} {J_{10}^6 }:   \frac { J_6^{10}} {J_{10}^6 } : 1
\right]
\end{equation}
So the triple 
\begin{equation}\label{abs-inv}
i_1=\frac {J_2^{30}} {J_{10}^6 },  \qquad i_2=  \frac {J_4^{15}} {J_{10}^6 },  \qquad i_3=   \frac { J_6^{10}} {J_{10}^6 }
\end{equation}
 uniquely determines the equivalence class of a genus 2 curve over the algebraic closure of $k$.  
 
Our initial intention was to  create a dictionary with keys $(i_1, i_2, i_3)$,  but these numbers blow up very quickly which makes and significant computations impossible. If the rational numbers have a significant number of decimal places, their exact representation might be lost when converted to floating-point format.
There is another set of invariants defined by Igusa in \cite{Ig},
\begin{equation}\label{t-invariants}
 t_1= \frac {J_2^5} {J_{10}},\quad  t_2 = \frac {J_4^5} {J_{10}^2},  \quad t_3 = \frac {J_6^5} {J_{10}^3}
\end{equation} 
which are defined everywhere in the moduli space.    They were used in some computations in \cite{arith-gen-2}, and are lower degrees than those in \cref{abs-inv}.   For relation among different kind of invariants  see \cite{m-sh, genus-2-univ} among many other places. 

We will use these invariants instead. Our Python dictionary will be keyed on ordered triples  $(t_1, t_2, t_3)$, whcich are  $\GL_2 (\bar k)$-invariants, hence every entry in the dictionary corresponds to the unique isomorphism class of genus 2 curves defined over $\bar \Q$. 
  
 
\begin{tabular}{c|c|c|c | c|}
 	& Entry 			& Value 					& Type 		& Description \\
 \hline
 0 	& $(x, y, z)$   		& $(t_1, t_2, t_3)$         		& float32  		& absolute  invariants        			\\  
 1 	& $\p$   			& $[J_2, J_4, J_6, J_{10}]$ 	& int      		& normalized moduli point			\\  
2 	& $\bar{\p}$   		&  						& int      		& absolutely normalized  point			\\   
3 	& wh     			& $\wh (\p)$ 				&   float32       	& weighted height  		\\  
4 	& awh     			& $\awh (\p)$                 		&   float32       	& absolute weighted height  		\\  
5 	& gcd     			& $\gcd (\p)$                 		&   float32       	& gcd of $\p$ 		\\  
6 	&  label1         		& T/F                        			& Boolean  	& True=fine, False=coarse \\
7 	& $[m, n]$			& 	$\Aut (\p)$				& [int, int] 		&  Gap Identity   			\\
8 	& label2  			& $\p \in {\mathcal L}_3$  				& Boolean   	&  \\
9 	& label3  			& $\p \in {\mathcal L}_5$  				& Boolean  	 &  \\
10 	& label4  			& $\p \in {\mathcal L}_7$  				& Boolean   	&  \\
\end{tabular}
 
 We  will describe later how to normalize points $\p=[a,b,c,d]$ in $\wP_{(2,4,6,20)}$.  
 
\subsection{Genus two curves with extra automorphisms}
The set of genus two curves  with extra automorphisms is a 2-dimensional irreducible subvariety of the moduli space which we will denote it by $\mathcal L_2$. It is the set of points   $\p \in \wP_{(2,4,6,10)}$  satisfying 

\begin{Small}
\[
\begin{split}
J_{30} (\p)=  	
&	41472 J_{10} J_4^5+159 J_4^6 J_2^3-236196 J_{10}^2 J_2^5-80 J_4^7 J_2+104976000 J_{10}^2 J_2^2 J_6-1728 J_4^5 J_2^2 J_6 \\
&	+6048 J_4^4 J_2 J_6^2  -9331200 J_{10} J_4^2 J_6^2-2099520000 J_{10}^2 J_4 J_6+12 J_2^6 J_4^3 J_6-54 J_2^5 J_4^2 J_6^2\\
& +108 J_2^4 J_4 J_6^3 +1332 J_2^4 J_4^4 J_6  -8910 J_2^3 J_4^3 J_6^2+29376 J_2^2 J_4^2 J_6^3-47952 J_2 J_4 J_6^4-J_2^7 J_4^4\\
& -81 J_2^3 J_6^4-78 J_2^5 J_4^5+384 J_4^6 J_6-6912 J_4^3 J_6^3	+507384000 J_{10}^2 J_4^2 J_2-19245600 J_{10}^2 J_4 J_2^3 \\
& -592272 J_{10} J_4^4 J_2^2+77436 J_{10} J_4^3 J_2^4+4743360 J_{10} J_4^3 J_2 J_6 -870912 J_{10} J_4^2 J_2^3 J_6 \\
& +3090960 J_{10} J_4 J_2^2 J_6^2-5832 J_{10} J_2^5 J_4 J_6-125971200000 J_{10}^3+31104 J_6^5+972 J_{10} J_2^6 J_4^2\\
&	+8748 J_{10} J_2^4 J_6^2-3499200 J_{10} J_2 J_6^3   = 0
\end{split}
\]
\end{Small}	
The locus $\mathcal L_2$ has two 1-dimensional loci corresponding to points with automorphism group $D_4$ and $D_6$; for details see \cite{deg2, m-sh, genus-2-univ, sha-2001,  sha-2016}.

\begin{figure}[h] 
   \centering
\includegraphics[scale=0.3]{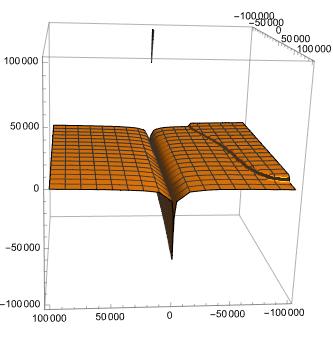}. \includegraphics[scale=0.3]{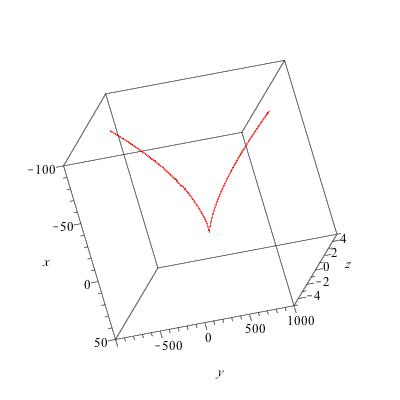}       \includegraphics[scale=0.3]{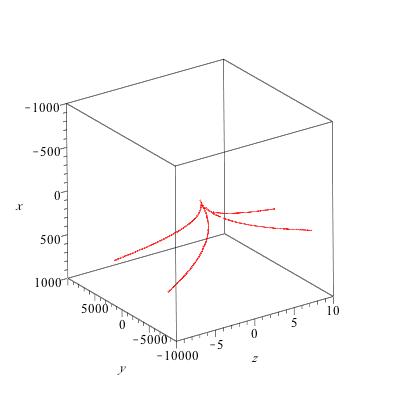}
   \caption{${\mathcal L}_2$ surface and 1-dimensional subloci of ${\mathcal L}_2$: curves with automorphism group $D_4$ or $D_6$}
   \label{fig:L2}
\end{figure}

\subsection{Genus two curves with $(n, n)$-split Jacobians}
Genus two curves with $(n, n)$-split Jacobians have been studied thoroughly during the last two decades and have been getting again some attention lately due to their use in isogeny based cryptography; see \cite{sha-2001, thesis, deg-3,   deg5,   sha-2016}
among others.   This extends the database from \cite{rat}. 

As in \cite{thesis}, we denote the locus of genus 2 curves with $(n, n)$-split Jacobian by ${\mathcal L}_n$. For $n$ odd, it is a 2-dimensional irreducible locus in $\M_2$.   Here we will describe how to create a database of points in ${\mathcal L}_n$ for $n=2, 3, 5, 7$.  

\begin{figure}[htbp] 
   \centering
   \includegraphics[width=2.6in]{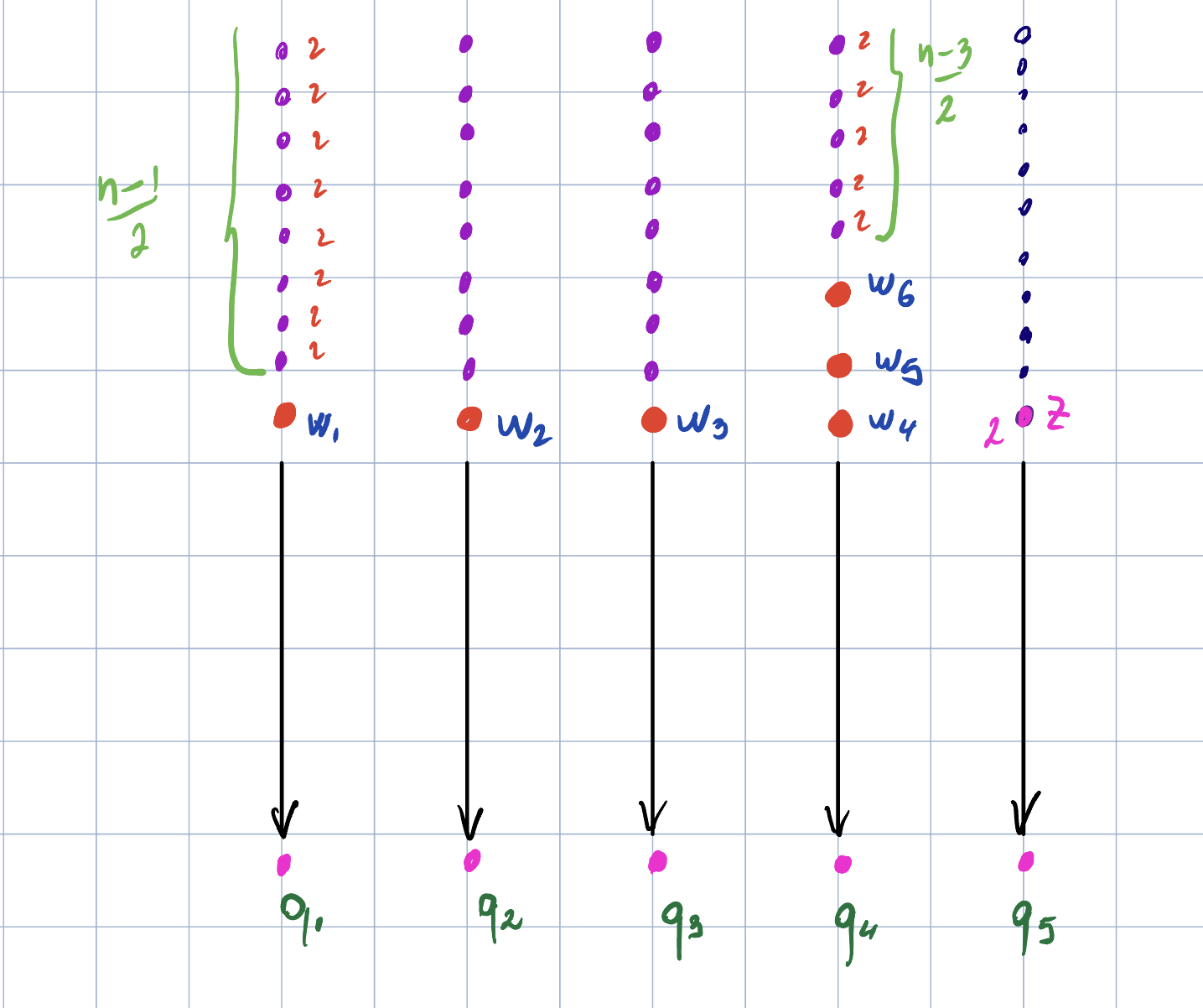} 
   \caption{Ramification type and Weierstrass points}
   \label{fig:ram}
\end{figure}

A degree $n$ covering $C\to E$, where $C$ is a genus two curve and $E$ an elliptic curve, induces a degree $n$ covering $\phi: {\mathbb P}^1 \to {\mathbb P}^1$ with ramification 
\[
\left(   2^{\frac {n-1} 2},  2^{\frac {n-1} 2},  2^{\frac {n-1} 2},  2^{\frac {n-3} 2},  2   \right)
\]
The unramified points in the fibers of the first four branch points are the Weierstrass points of the genus 2 curve; see \cref{fig:ram}.

Denote by $F_i (x)$, $i=1, \ldots 4$ the polynomial over the branch point $q_i$, which has as roots points of ramification index 2.  
Hence, $\deg F_1=\deg F_2= \deg F_2= \frac {n-1} 2$, and $\deg F_4= \frac {n-3} 2$. 
We fix a coordinate on the lower ${\mathbb P}^1$ by letting $q_1=0$, $q_2=\infty$, and $q_3=1$ and on the upper ${\mathbb P}^1$ by $w_1=0$, $w_1=1$, and $w_3=\infty$. Then
\[
\phi(x) = x  \left(    \frac { F_1 (x) }  {F_2 (x) }    \right)^2,   \quad \phi(x) - 1 = (x-1) \,  \left(    \frac { F_3 (x) }  {F_2 (x) }    \right)^2.
\]
Cases when $n=3$ and $n=5$ are special cases. When $n=3$ the fibers of $q_1$, $q_2$, $q_3$, and $q_5$ are identical, so we have extra symmetries permuting these branch points.   When $n=5$ then the fibers of $q_4$ and $q_5$ are identical so we have an extra involution permuting $q_4$ and $q_5$.   The first general case is when $n=7$ which is somewhat simpler since we don't have to worry about such extra symmetries, but of course everything becomes computational more challenging when $n$ gets bigger.  Next we provide formulas how to generate rational points in the moduli space of genus two curves, such that Jacobians are $(n,n)$-split for $n=2, 3, 5, 7$.  For justification and proofs of such formulas one can check \cite{thesis, deg2, deg-3, deg5}. 

\subsubsection{$(3,3)$-split:}   This case was studied in \cite{thesis} and summarized in \cite{deg-3}. 
Let $C$ be a genus 2 curve with $(3,3)$ split Jacobian.  Then from \cite[Theorem 4]{sha-2016} $C$ has equation
\begin{equation}\label{L3-param}
y^2= \left(v^2 x^3+u v x^2 + v  x +1\right)  \left(4 v^2 x^3 + v^2 x^2 + 2v x +1\right)
\end{equation}
for $\Delta= v\left( v -27\right) \left(4 u^{3}-u^{2} v -18 u v +4 v^{2}+27 v \right)  \neq 0$. 
Its moduli point is 
\[
\begin{split}
\p = &  \left[  2 v^4  \alpha  , 4 v^7  \beta , 4 v^{10}   \gamma, -16 v^{17} (v-27) \delta^3 \right] 
= 
\left[       2 v  \alpha  , 4 v \beta , 4 v   \gamma, -16 v^2 (v-27) \delta^3 \right]
\end{split}
\] 
since $v \neq 0$ and $\alpha, \beta, \gamma, \delta$ are 
\begin{Small}
\[
\begin{split}
\alpha  = & \,	4 u^{2}-12 u v +3 v^{2}+252 u -54 v -405  \\
\beta =  &	\,	u^{4} v -24 u^{4}-66 u^{3} v +9 u^{2} v^{2}+1188 u^{3}+297 u^{2} v +138 u \,v^{2}-36 v^{3}-8424 u v \\
& +945 v^{2}+14580 v \\ 
\gamma  =  &\,	2 u^{6} v^{2}-8 u^{5} v^{3}+2 u^{4} v^{4}-40 u^{6} v +106 u^{5} v^{2}+495 u^{4} v^{3}-204 u^{3} v^{4}+18 u^{2} v^{5}-144 u^{6} \\
& +1476 u^{5} v -18756 u^{4} v^{2}+4280 u^{3} v^{3}-1038 u^{2} v^{4}+564 u \,v^{5}-72 v^{6}+160704 u^{4} v \\
& +4464 u^{3} v^{2}+75024 u^{2} v^{3}-33480 u \,v^{4}+3186 v^{5}-104004 u^{3} v -1353996 u^{2} v^{2} \\
& +315252 u \,v^{3} -4032 v^{4}+3669786 u \,v^{2}-622323 v^{3}-2821230 v^{2}   \\
\delta  =   &\,	 4 u^{3}-u^{2} v -18 u v +4 v^{2}+27 v 
\end{split}
\]
\end{Small}
Notice that for rational values of $u, v$ we get rational points $\p \in {\mathcal L}_3$.  This provides an easy way to generate a database of  fine points in ${\mathcal L}_3$.
%
%
However, if we simply want   general rational points in ${\mathcal L}_3$, not necessarily fine moduli points, we have to use parameters $r_1, r_2$ which provide a birational parametrization of ${\mathcal L}_3$, but not necessarily genus 2 curves  defined over $\Q$; see \cite{deg-3} for details. 

In \cite{sha-2016} are given other examples of genus 2 curves with many elliptic subcovers  (Jacobian splits in more than one way). 
Recently curves with $(3, 3)$-split have been suggested by Flynn for genus two isogeny based cryptography; see \cite{flynn}.

\subsubsection{$(5, 5)$-split:} This case was studied in detail in \cite{deg5}. 
The parametric family of genus two curves in the locus ${\mathcal L}_5$ is given by
\begin{equation}\label{curve}
C: \qquad y^2 \ = \  x (x-1) \ \left(  a_3 x^3 + a_2 x^2 + a_1 x + a_0 \right) 
\end{equation}
where   
\begin{Small}
\[
\begin{split}
a_0 = \, & -b^4 (2 b^3 a+4 b^3-2 z a b^2+7 b^2 a^2+8 z b^2+4 b^2 +16 a
b^2+16 z b a+6 a^3 b+8 b a\\
& +2 z a^2 b+12 z b+16 b a^2+13 z a^2+z a^4+6 z a^3+4 z+12 z a)\\
a_1 = \, & -b^2 (12 b^3+12 b^4 a+32 z b a-6 a^4 b^2+44 b^2 a^3+6 b a^2+24 a
b^2+10 a^3 b+44 b^3 a^2+2 b a \\
& +52 b^3 a+61 b^2 a^2-12 b a^5-7 z a^2-2 z a+12 z b-4 a^6+12 b^4-a^4-40 z
a^3 b^2-16 z b^3 a^2 \\
& -12z a^5+36 z b^2-18 z a^3-26 z a^4+56 z a b^2+4 a z b^3+2 z a^2 b^2-20 z
a^3 b +28 z a^2 b\\
& +2 z a^6+24 z b^3+4 z b a^5-4 a^5-32 z a^4 b)\\
a_2 = \, &  5b^2 a^6+20 b^2 a^5+8 b a^6-61 b^4 a^2-18 b^5 a-56 b^4 a+4 z b
a+5 a^4 b^2-18 b^2 a^3 -24 z b^4\\
& -14z b^4 a-4 a b^2+8 b^3 a^4+2 b^3 a^5-54 b^3 a^3-70 b^3 a^2-24 b^3 a-14
b^2 a^2+4 a^4 b+10 b a^5 \\
& -6 z a^7+64 z a^3 b^3+38 z a^4 b^2+54 z a^3 b^2+12 z b^3 a^2-14 z a^6 b-10
z b^2 a^5-4 z a^7 b-4 a^6 z b^2 \\
& +32a^2 b^4 z+2 a^7 b-z a^8-36 z b^3-12 z a^5-12 z b^2-4 z a^4-28 z a
b^2-64 a z b^3-5 z a^2 b^2 \\
& +16 z a^2 b+28 z a^4 b-4 z b a^5-13 z a^6-12 b^5-12 b^4+34 z a^3 b\\
a_3  =\, & (2 a+1) (z a^4-2 a^3 b+4 z a^3+6 z a^3 b-4 b a^2+12 z a^2 b^2+10
z a^2 b-9 b^2 a^2+5 z a^2 \\
& -2 b a+2 z a-8 a b^2-12 b^3 a+8 a z b^3-4 b^3-4 z b-4 b^4-12 z b^2-8 z
b^3)\\
\end{split}
\]
\end{Small}
As mentioned above, there is an involution permuting branch points $q_4$ and $q_5$.  
Moreover,  $a, b$, and $z$ satisfy the equation
\begin{equation}\label{deg5-abz}
f(a,b,z):= \left(1+2 a \right) z^{2}+\left(-a^{2}-2 a b -2 a +2 b \right) z +2 a b +b^{2} =0 
\end{equation}
see  \cite[Thm. 2]{deg5} for details. 

\smallskip

\noindent \textbf{Equation of ${\mathcal L}_5$:} 
Computing the equation of ${\mathcal L}_5$ sounds as an easy exercise in elimination theory:  compute $i_1 (t, s), i_2(t,s), i_3(t,s)$ and eliminate $t$ and $s$.  This will give an affine equation $f(i_1, i_2, i_3)=0$. To get the projective equation, substitute $i_1, i_2, i_3$ in terms of $J_2, J_4, J_6, J_{10}$. 

The main problem with the above approach is that the degree of rational functions $i_1, i_2, i_3$ in terms of $t$ and $s$ are very large, which makes the elimination of $t$ and $s$ practically impossible.  

In \cite{deg5} was given the following approach in the computation of ${\mathcal L}_5$.  Let 
\[
u = \frac {2 a \ (a b+b^2+b+a+1)} {b \ (a+b+1)},  \qquad v= \frac {a^3} { b \ (a+b+1)}, \quad  w = \frac {(z^2-z+1)^3} {z^2 (z-1)^2}
\]
They are invariants of a group action on $k(a,b,z)$.  Since the modular invariants $J_2, J_4, J_6, J_{10}$ are invariants of any permutation of fibers, they can be expressed in terms of $u, v, w$.  Moreover,
\[k ({\mathcal L}_5) = k (u, v, w), \]
where the equation of $w$ in terms of $u, v$ is
\begin{equation}\label{u-v-w}
 c_2 w^2 + c_1 w + c_0 =0
\end{equation}
with $c_0, c_1, c_2$ as follows:
\begin{equation}\label{eq_w}
\begin{split}
c_2  = & \,  64v^2(u-4v+1)^2\\
c_1  = & \, -4v(-272v^2u-20 v u^2+2592 v^3-4672 v^2+4 u^3+16 v^3 u^2-15 v u^4\\
&\,  -96 v^2 u^2+24 v^2 u^3+2 u^5-12 u^4+92 v u^3+576 v u-128 v^4-288 v^3 u)\\
c_0  = & \, (u^2+4 v u+4 v^2-48 v)^3 \\
\end{split}
\end{equation}
Notice that the surface in $u, v, w$ is a septic surface and more difficult to get a parametrization of it. However, expressing $i_1(u,v,w), i_2(u,v,w), i_3(u,v,w)$ together with the \eqref{u-v-w} gives us a system of equations which is easier to handle and possible to eliminate $u, v, $ and $w$.   In \cite{deg5} was shown that   $k(u, v, w)= k ({\mathcal L}_5)$.  

\begin{rem}
The equation of ${\mathcal L}_5$ was computed in \cite{deg5} by the above method and using resultants to eliminate $u$ and $v$. Since it is too long it was not displayed in the paper, but in a webpage that no longer is available.  Due to regular and repeated requests for this equation, we intend to display it on our groups's webpage, even though it is exactly  its impracticality  of usage due to its length which motivates this paper. 
\end{rem}

\subsubsection{Cases for $n\geq7$}    Both cases $n=3$ and $n=5$ are special cases due to their ramification structure.  For example, for $n=3$ the fibers $\phi^{-1} (q_i)$, for $i=1,2,3$ and $\phi^{-1} (q_5)$ are the same and for $n=5$ the fibers  $\phi^{-1} (q_4)$ and $\phi^{-1} (q_5)$ are the same.  This fact induces extra symmetries and therefore a group action as shown in the computation of these spaces.  

The first case which is a general case (i.e. the ramification structure is the same as large $n$)   is the case $n=7$.  Its computation is much more involved; see \cite{sha-2001} for computations of some of it degenerate loci.     There are other ways how to generate rational points in ${\mathcal L}_n$, even though not general points.  The general ramification for $n>7$ has four cases (called degenerate cases in \cite{deg-3, deg5}).  These cases give 4  curves in ${\mathcal L}_n$ and two of these curves are genus zero curves.  Hence, a rational parametrization of these curves would provide rational points in ${\mathcal L}_n$. However, any model based only on these points would be suited only for this curve and not the whole ${\mathcal L}_n$ space.   In \cite{sha-2016, b-e-sh} was considered the cases of intersection between different $\mathcal L_n$
In general, computing the $\mathcal L_n$ locus requires Groebner bases and such methods are inefficient as $n$ increases, which makes the machine learning methods even more appealing.


\section{Number of points in ${\mathbb P}_\w (\Q)$ with bounded weighted height} 
In this section we want to estimate the number of points in a weighted projective space. Since our main application is the moduli space of hyperelliptic curves, we focus on the binary forms. 
Let  ${\mathcal R} _d$  be the ring of invariants of degree $d\geq 3$ binary forms and $\xi=(\xi_0, \ldots , \xi_n)$ be a generating  of ${\mathcal R} _d$    with  weights  $\w=(q_0, \ldots , q_n)$. 
Let  $ {\mathbb P}_\w^n$  denote the weighted projective space over $\C$  (say ${\mathbb P}_{\w, \C}^n$) and ${\mathbb P}_w^n (\Q)$ the set of points with rational coordinates.  Notice that for each $\p \in {\mathbb P}_\w^n (\Q)$ we can assume $\p = [x_0, \ldots , x_n ]$, where $x_i \in \Z$, for all $i=0, \ldots , n$. 

Fix $h\in \R^{\geq 0}$ and let 
\[
\begin{split}
B_h & := \{ \p \in {\mathbb P}_\w^n (\Q)   \; \mid \; \awh ( \p ) \leq h) \}, \\
C_h  & := \{ \p \in {\mathbb P}_\w^n   (\Q)    \; \mid \; h-1 < \awh ( \p ) \leq h) \}, 
\end{split}
\]
Recall that the wighted height is defined for all $\p \in {\mathbb P}_\w^n (\Q)$ by taking a normalized representative of  $\p$ in   ${\mathbb P}_\w^n (\Z)$ and then applying definition in \cref{def:height}.

Let $F_d : [0, \infty) \to \Z^{\geq 0}$ be the function which denotes the cardinality of $B_h$.  In other words, $F_d (h)$ is the number of equivalence classes of degree $d\geq 2$ binary forms with absolute weighted height less or equal to $h$.  
Similarly, let    $G_d : [1, \infty) \to \Z^{\geq 0}$ be the cardinality of $C_h$ or equivalently the number of equivalence classes of degree $d\geq 2$ binary forms with absolute weighted height  in $(h-1, h]$. 
Hence, 
\[
G_d(h) = F_d (h) - F_d (h-1).
\]
From   \cite[Theorem 1]{b-g-sh} we know that    $F_d(h)$ and $G_d(h) $ are  well defined. 

\begin{thm}\label{thm1}   
Let $\w=(q_0, \ldots , q_n$ and  $\wP_{\w, \Q}^n$ a well-formed weighted projective space. 
The number of points in   $\wP_{\w, \Q}^n$ with height less or equal to a number $h$ is   
\[
F_d (h) \leq  \sum_{i=0}^n   \left(  h^{q_{n-i}} \cdot \prod_{j=0}^{n-i} \left( 2 h^{q_j}+1 \right)    \right)
\]
\end{thm}

\proof
Assume $\awh (\p) \leq h$.    For each $i=0, \ldots n$ we have 
\[
|x_i|^{\frac 1 {q_i}} \leq h \implies |x_i| \leq  h^{q_i}
\]
Hence, there are $2 h^{q_i} + 1$ choices for $x_i$. Moreover, we can normalize one of the coordinates so it is always positive. We can do that for the highest power, which by our ordering     is $q_n$. Then there will be only $h^{q_n}+1$ choices for that coordinate. 
Hence, our total number is bounded by 
\[
\left(    h^{q_n} +1 \right) \cdot \prod_{i=0}^{n-1} \left( 2 h^{q_i}+1 \right)
\]
Since we are counting stable points in the moduli space of binary forms, then at least one of the coordinates must be nonzero; see \cite{curri}.  Assume $x_n\neq 0$.  Then there are $h^{q_n}$ choices for $x_n$ and the number of such forms is less than 
\[
  h^{q_n} \cdot \prod_{j=0}^{n-1} \left( 2 h^{q_j}+1 \right).  
\]
If $x_n=0$ then by the same argument  there are 
\[
 h^{q_{n-1}} \cdot \prod_{j=0}^{n-2} \left( 2 h^{q_j}+1 \right).
\]
and so on. Adding up all the cases we have
\[
F_d (h) \leq  \sum_{i=0}^n   \left(  h^{q_{n-i}} \cdot \prod_{j=0}^{n-i} \left( 2 h^{q_j}+1 \right)    \right)
\]
This completes the proof.
\qed

Notice that the above method counts as different moduli points tuples $\lambda \star (x_0, \ldots , x_n)$ as long as $\max \{ |\lambda^{q_i} x_i|^{1/q_i}\} \leq h$. Formula is precise only for $h=1$ as it will be seen in computations with binary sextics. Consider for example 
$[1:0:0:0]$ and $[2: 0:0:0] $ in ${\mathbb P}_{(1,2,3,5)}$.  They are the same point but counted separately  by the above formula when $h > \sqrt{2}$. 

Next we use the results of the above two sections to study the moduli space $\M_2$ of genus 2 curves.  $\M_2$ has been the focus of investigation for a long time due to its particular place in the algebraic geometry of curves and its applications to hyperelliptic curve cryptography.  In particular we want to compare our result to those in \cite{rat, m-sh}.   
We will continue to use the terminology as above. 

Let ${\mathcal H}_g$ denote the moduli space of genus $g \geq 2$ hyperelliptic curves defined over $k$ and $B_{2g+2}$ the space of degree $2g+2$ binary forms (up to equivalence). 
  Since every genus $g$ hyperelliptic curve corresponds to a binary form of degree $2g+2$ we have an embedding
\[
\phi:    {\mathcal H}_g \hookrightarrow B_{2g+2}
\]
The image $\phi  ({\mathcal H}_g)$ does not intersect with the locus $\Delta=0$, where $\Delta$ is the locus of all binary forms with zero discriminant.  So ${\mathcal H}_g$ is isomorphic to $\wP_{\w}^n \setminus \{\Delta =0\}$, where $\w$ is the set of weights corresponding to the degrees of the generators of the ring of invariants ${\mathcal R}_{2g+2}$.   Hence, to count points in the moduli space ${\mathcal H}_g$ it is enough to count equivalence classes of degree $(2g+2)$ binary forms with nonzero discriminant. 

\begin{lem}
The number of points in ${\mathcal H}_g$ with weighted moduli height $\leq h$ 
\[
h^{q_n} \cdot    \sum_{i=0}^{n-1}   \left(  h^{q_{n-i}} \cdot \prod_{j=0}^{n-i} \left( 2 h^{q_j}+1 \right)    \right)
\]
\end{lem}

\proof
The invariant with the highest degree is the discriminant of the binary form.  Since this binary form correspond to a hyperelliptic curve, the discriminant is not zero. Hence, $x_{q_n}\neq 0$.  Then, there are only $h^{q_n}$ choices for $x_{q_n}$; see the proof of \cref{thm1}. This completes the proof. 
\qed

\begin{table}[h]
\centering
\caption{Points   $B_6\setminus\{ \Delta=0\}$ or equivalently the number of genus 2 curves with moduli height $\leq h$.}
\label{table-x}
\begin{tabular}{c|c}
\toprule
 h & \# of points in $B_6$ \\
\midrule
1 & 40 \\
2 & 24 862 \\
3 & 1 781 202 \\
4 & 39 251 668 \\
5 & 440 104 780 \\
6 & 3 195 496 050 \\
7 & 17 146 927 462 \\
8 & 73 657 853 512 \\
9 & 266 816 523 888 \\
10 & 844 626 323 110 \\
\bottomrule
\end{tabular}
\end{table}

 Next we focus on our special case when $g=2$.  In this case, as mentioned above, the wieghts are $\w=(2,4,6,10)$ and the invariants $J_2$, $J_4$, $J_6$, and $J_{10}$.  We can count the points in  ${\mathbb P}_w (\Q)$ using the above approach and having in mind that $J_{10}\neq 0$.  

\begin{lem}
The number of points in $\wP_{(1,2,3,5)}$ with weighted moduli height $\leq h$ 
\[
\begin{split}
F_6(h) & \leq   h^{5}  {\left(2 \, h^{3} + 1\right)} {\left(2 \, h^{2} + 1\right)} {\left(2 \, h + 1\right)}  +
h^{3}   {\left(2 \, h^{2} + 1\right)} {\left(2 \, h + 1\right)} 
 +h^{2}  {\left(2 \, h + 1\right)}  + h \\
& = h {\left(8 \, h^{10} + 4 \, h^{9} + 4 \, h^{8} + 6 \, h^{7} + 2 \, h^{6} + 6 \, h^{5} + 3 \, h^{4} + 2 \, h^{3} + 3 \, h^{2} + h + 1\right)} 
\end{split}
\]
Moreover, there are exactly 27 genus 2 curves with weighted moduli height $\awh =1$.
\end{lem}

\proof
The invariant with the highest degree is the discriminant of the binary form.  Since this binary form correspond to a hyperelliptic curve, the discriminant is not zero. Hence, $x_{q_n}\neq 0$.  Then, there are only $h^{q_n}$ choices for $x_{q_n}$; see the proof of \cref{thm1}. This completes the proof. 
\qed

In   Table~\ref{table-0}  we display  all points $\p \in  \mathbb{WP}_{(1,2,3,5)}^3 (\Q) \setminus  \{J_{10} = 0\} $ with weighted moduli heights $\h = 1$.  

\begin{small}
\begin{table}[htb]
\centering
\caption{Moduli points $\p=[J_2: J_4: J_6: J_{10}]$ with  height $\awh=1$}
\begin{tabular}{|c|c|c|c|c|c|}
	\toprule
	\# & $\p$    & \#   &    $\p$  & \#   &    $\p$  \\
	\midrule
	1 & [0, -1, 0, 1]  & 	10 &  [1, 0, 1, 1]  &	19 & [0, 1, 1, 1] \\
	2 & [0, 1, 0, 1]    & 	11 & [1, -1, -1, 1] & 	20 & [1, 0, 1, -1]\\
	3 & [0, -1, 1, 1] & 12 &  [1, 1, -1, 1] & 21 & [1, -1, -1, -1]  \\
	4 & [0, 0, 0, 1] & 13 &  [1, 1, 1, -1] & 22 & [1, 1, -1, -1]\\
	5 & [0, 0, 1, -1] & 14 & [1, -1, 1, -1] & 23 & [1, -1, 0, -1]\\
	6 & [0, 0, 1, 1] & 15 & [1, 1, 1, 1] & 24 & [1, 1, 0, -1]\\
	7 & [1, 0, -1, 1]   & 16 & [1, 0, -1, -1] & 25 & [1, 1, 0, 1]\\
	8 & [1, 0, 0, -1]  & 17 & [0, -1, 1, -1] & 	26 & [1, -1, 0, 1] \\
	9 & [1, 0, 0, 1] & 18 & [0, 1, 1, -1] & 	27 & [1, -1, 1, 1] \\
	\bottomrule
\end{tabular}
\label{table-0}
\end{table}
\end{small}

It is worth noting that sometimes, for different reasons,  we use the weighted projective spaces which are not well formed, even though every weight projective space is isomorphic to a well formed one; see \cite{s-sh}.  However, these counting functions $F_d (h)$ and $G_d (h)$ do not take the same values on such spaces.  Consider for example the isomorphism in \cref{iso-wp}.   The weighted height of a point in $\wP_{(1, 2,3,5)}$ is the square of the weighted height of a point in $\wP_{(2,4,6,10)}$.  

\begin{lem}
In the case of binary sextics  the value of  $G_6 (h)$  in $\wP_{(1,2,3,5)}$ is 
\begin{equation}
\begin{split}
G_{   \wP_{(1,2,3,5)}  } (h)  =    \,	& 88 h^{10}-400 h^{9}+1176 h^{8}-2256 h^{7}+3038 h^{6}-2862 h^{5}\\
			& +1879 h^{4}-812 h^{3}+215 h^{2}-28 h +2
\end{split}
\end{equation}
Moreover,  for the Veronese embedding $\phi:  \wP_{ (2,4,6,10) }   \to \wP_{(1,2,3,5)}  $ as in \cref{iso-wp} we have 
\[
 \awh_{   \wP_{(2,4,6,10)}  }  (\p)   = \left(   \awh_{   \wP_{(1,2,3,5)}  }  (\phi (\p))   \right)^{\frac 1 2}
\]
which implies that 
\[
G_{   \wP_{(1,2,3,5)}  }  (1)   = G_{   \wP_{(2,4,6,10)}  }  (1)
\]
\end{lem}

\proof The proof is an imediate consequence of the above lemma, definition of the weighted height, and the isomorphism in  \cref{iso-wp}. \qed
 
\begin{cor}
The number of points in $\wP_{(2,4,6,10)} (\Q)$    with weighted height $\leq \awh$ is 
\[
F_{   \wP_{(2,4,6,10)}  }  (h)   =  F_{   \wP_{(1,2,3,5)}  }  (h^2)
\]
\end{cor}

Now we have a rough estimate on the number of points in $\wP_{(2,4,6,10)}$, which is our main focus of study. 

\section{A database of genus two curves}\label{sec-data}
 Now that we know how to generate rational points in ${\mathcal L}_n$, we would like to see if we can generate some random data in $\M_2$ and train a model that answer arithmetic properties of $\p \in \M_2$, including whether $\p \in {\mathcal L}_n$.    Our data will have points from ${\mathcal L}_n$, for $n=2,3,5,7$ and all points of weighted moduli height $\awh \leq 3$.   
  
\subsection{Creating the database}
An entry in the dictionary looks like:
\[
(x,y, z) :  (  \p, \awh (\p), \text{Fine},  \Aut (\p), \p \in {\mathcal L}_3, \p\in {\mathcal L}_5, \p \in {\mathcal L}_7 )
\]
where   
\[
\begin{split}
\p & = [J_2, J_4, J_6, J_{10} ]   \; \text{weighted moduli point} \\
\awh & = \textit{absolute weighted height} \\
\text{Fine} &= \textit{ True/False }\\
\Aut & = \text{Automorphism group} \, \Aut (\p) \\
flag & = (3, 3)-\text{split Jacobian} \\
flag & = (5, 5)-\text{split Jacobian} \\
flag & = (7, 7)-\text{split Jacobian} \\
\end{split}
\]
where the key $(x, y, z)$ is the triple of absolute invariants $(i_1, i_2, i_3)$.  Notice that the automorphism group basically determines if the corresponding point is in the ${\mathcal L}_2$ locus or not.  Other then the curve $y^2= x(x^5-1)$, all curves with $|   \Aut( \p) |   >2$ are in the ${\mathcal L}_2$ locus and they are fine moduli points. 

\subsubsection{All points with weighted height $\leq 2$}  First we compute all points $\p$ with weighted height $\awh (\p) \leq 2$.  From our estimates  in the previous section  there are $\leq 39 251 668$   points in $\wP_{(2,4,6,10)}$.   

\subsubsection{Generating points in the locus ${\mathcal L}_2$}
Unfortunately there aren't many points from ${\mathcal L}_2$ with weighted height $\leq 2$; see \cref{tab:L2}. However, we can generate many other rational points in ${\mathcal L}_2$ by using the birational $(u, v)$- parametrization in \cite{deg2}.

\subsubsection{Generating points in the locus ${\mathcal L}_3$}
We use the Eq.~(8) and for each $(u, v)$-value compute the moduli point in terms of Igusa invariants. Using this method we get fine moduli points.  If we want random rational points in ${\mathcal L}_3$ (not necessarily fine) we can use the birational $(r_1, r_2)$-parametrization as in \cite{deg-3}. 

\subsubsection{Generating points in the locus ${\mathcal L}_5$}
The  cubic surface in \cref{deg5-abz} is quadratic in $z$ and therefore     rational, by a result of Clebsch.  Hence $a, b \in k(t, s)$ for some parameters $t$ and $s$.
Thus invariants $J_2, J_4, J_6, J_{10} \in k(t, s)$. Giving random values to $t$ and $s$ generates rational points in ${\mathcal L}_5$. 

Notice that for all $z=s \neq 0, 1$, the curve $f(a, b, s)$ is a genus zero curve for which we can get a parametrization $\left( a(t), b(t) \right)$. Hence, we get the desired parametrization in $t, s$ for the surface.  

\begin{exa}
For $s=2$  we get $a=-\frac{8}{t^{2}+2 t -2}$ and $b=-\frac{2 \left(t^{2}-2 t -2\right)}{t^{2}+2 t -2}$  and now can randomly pick $t\in \Q$ such that the denumerator is $\neq 0$.

For $s=1/2$   we get $a =  \frac{16}{t^{2}-4 t -8}$ and   $b=-\frac{t^{2}+4 t -8}{2 \left(t^{2}-4 t -8\right)  }$. 
\end{exa}

We get the following algorithm.

\begin{algorithm}
\caption{Algorithm to generate $n$ rational points in $\p \in {\mathcal L}_5$ with $\awh (p)\leq h$.}
\begin{algorithmic}
\STATE \textbf{Input:}  Integer $n$.    
\STATE \textbf{Output:} A list of $n$ rational points in  ${\mathcal L}_5$
\STATE
\STATE    Generate a list of random rational numbers $S=\{ s_1, \ldots , s_r \}$, for $r =10$
\FOR{each number $i$ in $1, \dots , r$}
	\STATE Randomly pick a rational number $s\in \Q$
    \STATE Parametrize $f(a,b,s)=0$ in a parameter $t$
    \STATE Generate   $n/10$ rational points  $(s, a(t), b(t))$ 
    \STATE For each $(s, a(t), b(t))$  compute $\p_{s, t}= (J_2, J_4, J_6, J_{10})$
\ENDFOR
\STATE \textbf{return} The set of points  $\{ \p_{s, t} \}$
\end{algorithmic}
\end{algorithm}

\subsection{Normalizing the data}
Our data is given in terms of tuples of Igusa invariants.  
There are two main ways of normalizing the data using two functions \textsl{MinIgusaTuple ()} and \textsl{MinAbsIgusaTuple ()}. We briefly explain the differences between the two.  

Let $\p=[x_0 : x_1: x_2 :x_3]\in \wP_{(2,4,6,10)} (\Q)$. Without loss of generality we can assume that its coordinates $x_0, \ldots , x_3\in \Z$. A weighted moduli point is called \textbf{normalized} if    $\wgcd(x_0, \ldots , x_3)=1$ and \textbf{absolutely normalized} if 
$\awgcd(x_0, \ldots , x_3)=1$.

The function \textsl{MinIgusaTuple ($\p$)} converts any point $\p$ to a normalized point.  
The function \textsl{MinAbsIgusaTuple ($\p$)} converts any point $\p$ to an absolutely normalized point.    The algorithm to do this involves integer factorization which makes it  inefficient for our data with large coordinates $J_2, J_4, J_6, J_{10}$. 

There have ben attempts to find a polynomial time algorithm for determining the weighted greatest common divisor of a tuple. The following elementary result shows that this can't be done. 

\begin{lem}   Let $\x = (x_0, \ldots , x_n)\in \Z^{n+1}$.  Then  $\wgcd (\x)  | \gcd (x)$. 
Determining the weighted greatest common divisor of a tuple is at best equivalent to integer factorization. Equivalently, 
normalizing in a weighted projective space is equivalent to integer factorization.
\end{lem}

\proof   The first part comes from the definition of weighted greatest common divisor. For the second part, 
take the tuple $\x= (d, \ldots, d)$.  Then $\gcd (\x)=d$ and to find the $\wgcd (\x)$ we need to know all prime factors of $d$. 
\endproof

It is worth noticing that for most of our database, which includes points of weighted heights $\leq 5$ the gcd of tuples does not have prime factors or it has small primes as factors so the normalizing can be done rather quickly.

\subsection{Distribution of fine points in the moduli space}
There are two types of points in the weighted moduli space $\wP_{(2,4,6,10)}$, namely \emph{fine points} and \emph{coarse points}. Fine points are those points such that their field of moduli is a field of definition, while the rest of points are called coarse points. 

We also classify fine points in two classes, those with extra automorphisms and those which have automorphism group isomorphic to the cyclic group of order 2.   

Let $\x\in {\mathbb P}^n_{\w} (\Q)$.  Can we find a degree $d$ binary form $f \in B_d$, defined over $\Q$,  such that $\x= [f]$.  The answer is in general negative. Points for which we can answer positively the above question are called \textbf{fine moduli points}, otherwise $\x$ has a non-trivial obstruction and we will call it a \textbf{coarse point}.  The problem of determining which points are fine points is referred to as field of moduli versus field of definition problem in arithmetic of moduli spaces of curves; see \cite{m-sh} among many others. 
\begin{figure}[h!]
\begin{center}
\includegraphics[scale=0.35]{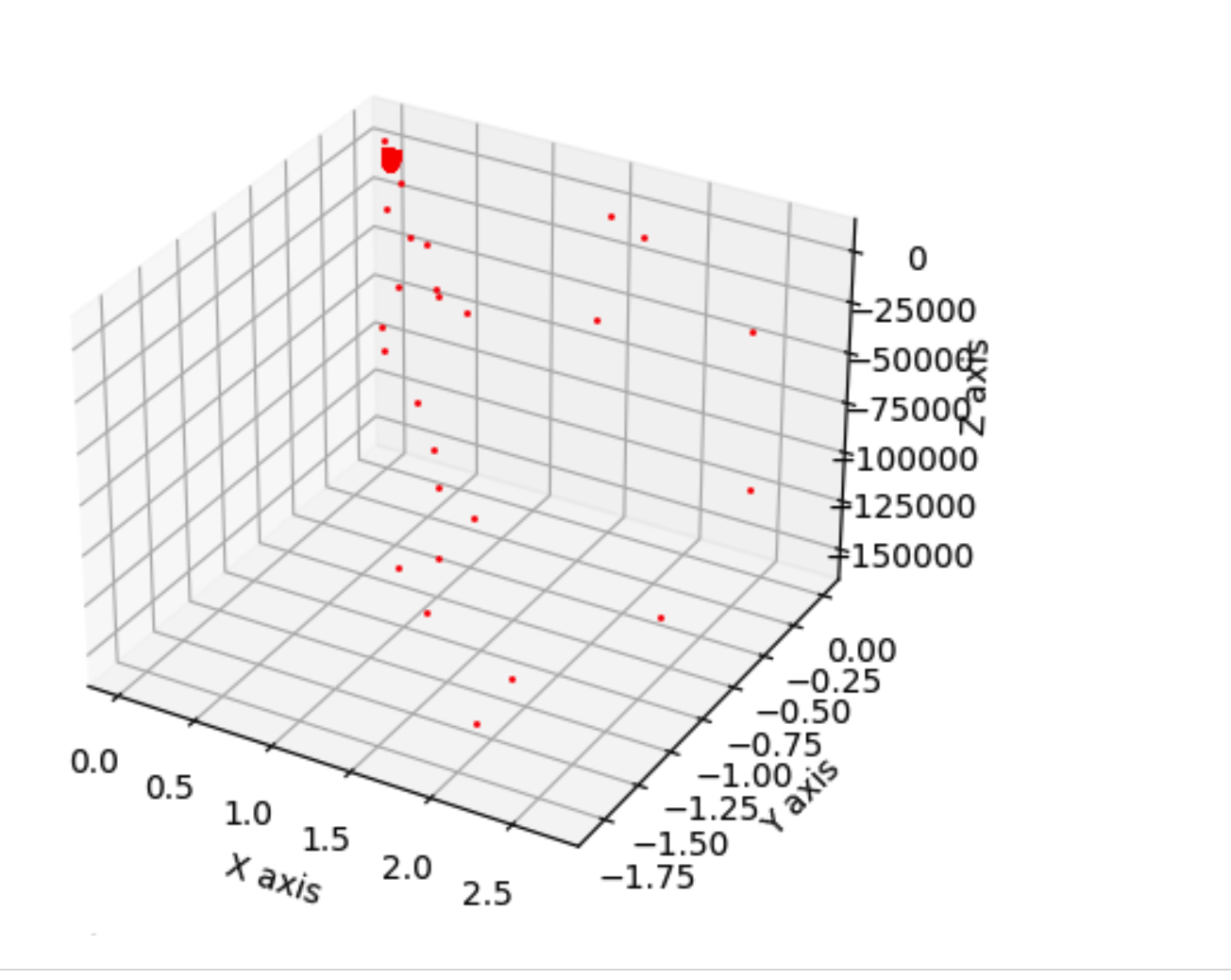}
\caption{Distribution of points with extra automorphisms}
\label{Aut}
\end{center}
\end{figure}
 Let us now see what is the distribution of fine points with automorphisms in our database (red points).  From a computational point of view it is quite hard to do this simply by brute force for our database which has about 500 000 points. Instead we will use the above models to see what information we can gather and then prove our results (if any) via brute force computationally. 
 By taking a random sample and graphing all red points we get the   picture \cref{Aut}.
Red points seem to be very scarce around the origin of the coordinate system and secondly there are no green points (fine points with trivial automorphism group).   Somewhat to be expected by people who have extensive computational experience with the moduli space of genus 2 curves, but not any obvious  reason for it.

\subsubsection{Coarse moduli points in $\M_2$}
Let us first go over some of the conditions that a point $\p\in \wP_{(2,4,6,10)}$ is a fine moduli point.    Let  $\p\in \wP_{(2,4,6,10)}$  be a rational point. Then  a genus 2 curve $\X$ is defined over $\Q$ if and only if the conic $Q= \X/\, \langle w\rangle$ has a rational point, where $w$ is the hyperelliptic involution. 
 From \cite[Lem. 3.1]{genus-2-univ} the conic $Q$ is isomorphic to the diagonal conic 
 \[
 Q^\prime : \; x_1^2 - \gamma x_2^2  - \Lambda_6 x_3^2 =0,
 \]
 where $\gamma$ and $\Lambda_6$ are determined in terms of Igusa invariants as in \cite{genus-2-univ}.  This conic has a rational point if and only if 
 there exist rational numbers $\alpha, \beta \in \Q$ such that   
 \[
 \alpha^2 + \left( \Lambda_6 \sigma\right) \beta^2 = \gamma.
 \]
 where $\Lambda_6, \sigma, \gamma$ are given explicitly in \cite{genus-2-univ} in terms of Siegel modular forms or equivalently in terms of the invariants $J_2, J_4, J_6, J_{10}$. 
A random point  $\p\in \wP_{(2,4,6,10)}$ generates random values for $\Lambda_6, \sigma, \gamma$, which
can not be written as a sum of squares of two rational numbers $\alpha$ and $\beta$. Hence we have the following. 

\begin{lem}
A generic point in $\M_2$ is a coarse point.   
\end{lem}

We graph below all   points for weighted height $\wh (\p) \leq 3$.    Blue points are overwhelming as expected  from the above Lemma. 
\begin{figure}[h] 
   \centering
   \includegraphics[width=1.6in]{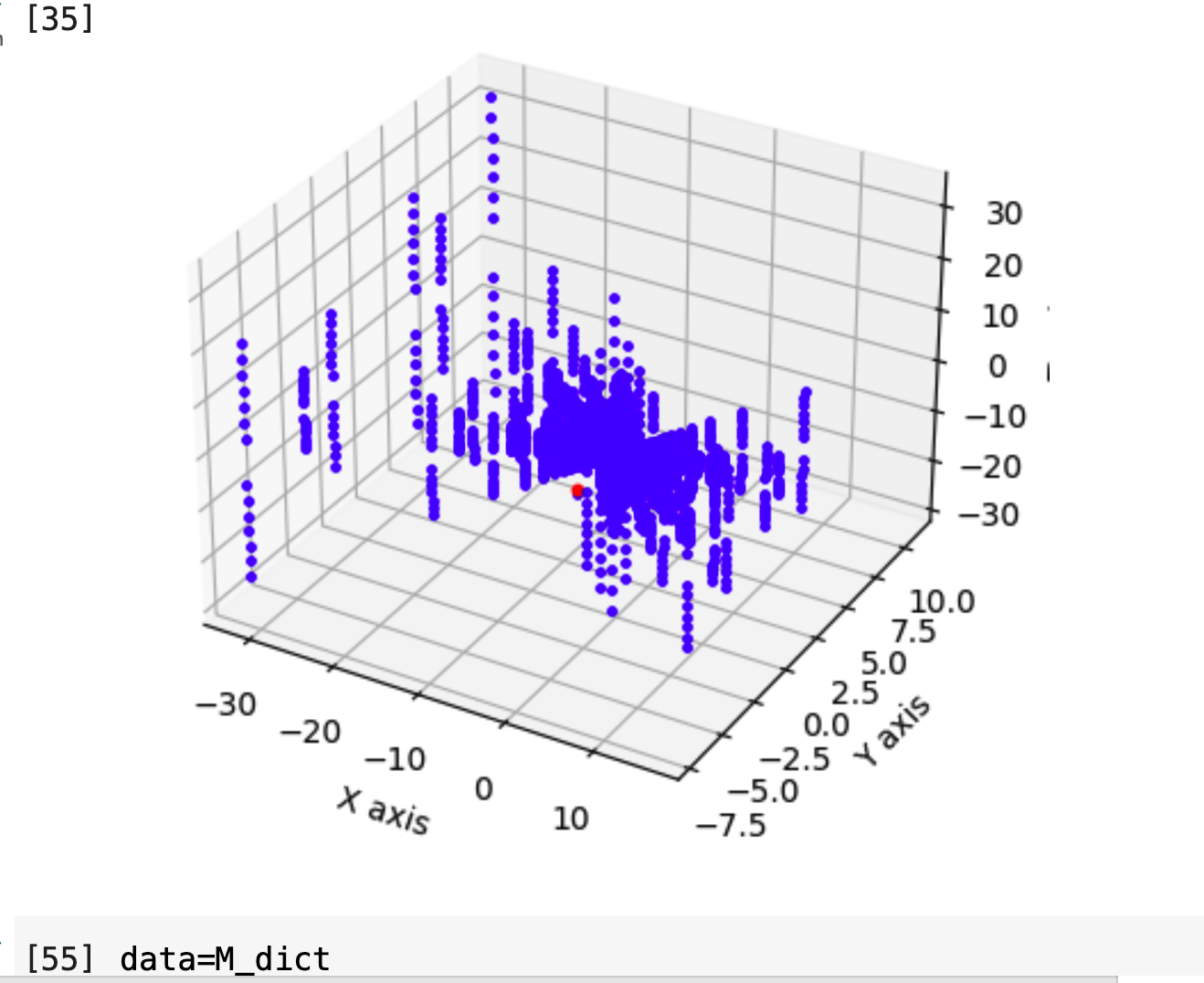} \includegraphics[width=1.6in]{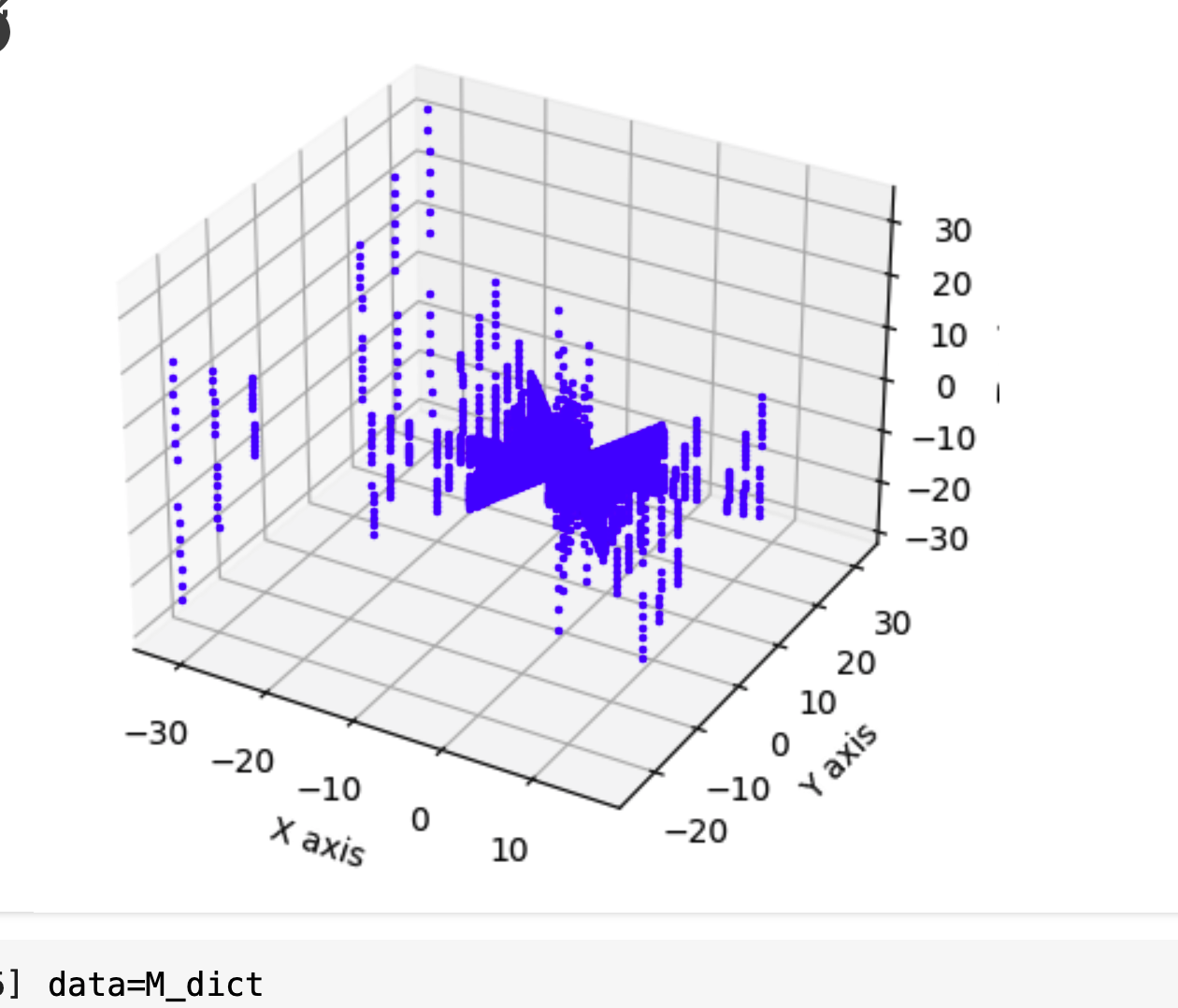} \includegraphics[width=1.6in]{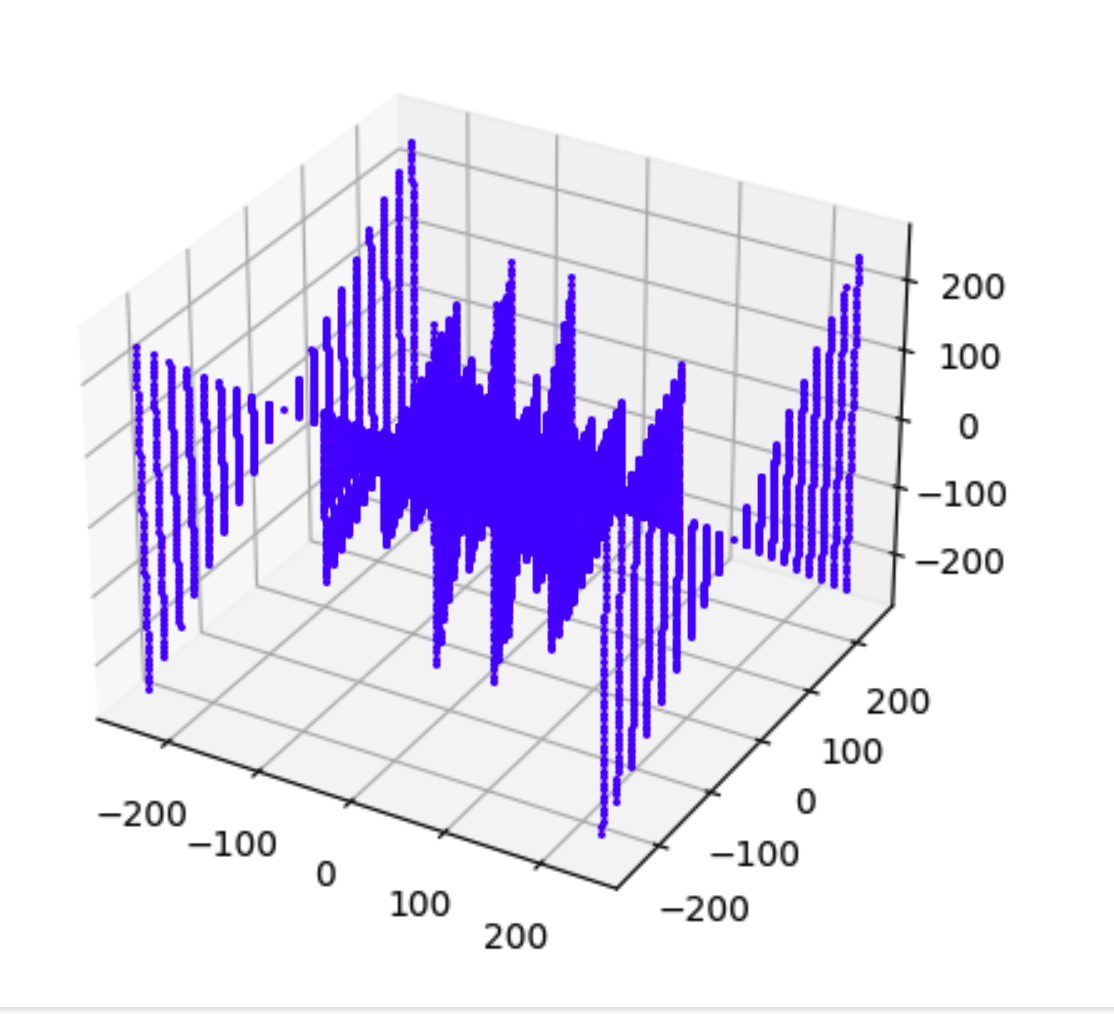} 
   \caption{Graph of rational points of weighted height $\wh \leq 2$}
   \label{h-1}
\end{figure}
The graphs above led us to check if our models are not working correctly or if this is truly the case.   

Perhaps  a larger database would make possible to  graphically see the distribution of red points (fine points with extra automorphisms) and green points (fine points with automorphism group of order 2).     Overall, it seems as for small weighted moduli height there are no such points.  

Surprisingly there is only one red dot in all these graphs. These graphs were obtained using the sequential method and the existing red dot was no surprise because that is a very special genus 2 curve and well known, namely the curve with $J_2=J_4=J_6=0$ which correspond to the single moduli point $\p$ with $\Aut (\p)\cong C_{10}$.   

Hence, if the graphs are accurate it would mean that there is no point of ${\mathcal L}_2$ in our database for weighted height $\leq 2$  (all points in ${\mathcal L}_2$ are red points).   Hence,  these graphs  were  a strong enough reason for us to go through all the cases for $\wh (p)\leq 3$ and check computationally if there are points in ${\mathcal L}_2$ or ${\mathcal L}_3$.   

\begin{rem}
Notice that this fact alone (scarcity of genus two curves in $\mathcal L_2$ or $\mathcal L_3$ with small Igusa invariants)  had never been noticed before and probably would have never been noticed if it wasn't for such approach.  It emphasizes our philosophy that machine learning techniques could be a great tool in detecting mathematical patterns which can then be proved with more classical methods. 
\end{rem}

\begin{lem}
There are no rational points $\mathfrak p\in {\mathcal L}_2$ with weighted  height $\wh (\p) < 3/2$. 
\end{lem}

\begin{table}[h!]
\caption{Rational points of height $\leq 3$ in ${\mathcal L}_2$}
\label{tab:L2}
\begin{tabular}{|c|c||c|c||c|c|}
\hline
 \# & $\p$   &   \# & $\p$ &   \# & $\p$  \\
\hline
1 & [4, -14, 2, 1] 	& 2& [2, -11, 5, 1]   	& 3 & [-2, -8, 14, 1] \\
4& [-2, 16, -14, 1] 	& 5& [2, 13, -3, 1] 	&  6& [4, 16, 0, 2] \\ 
7& [4, -8, 16, 2] 	& 8& [0, -3, 27, 2]   	& 9& [-4, 4, 28, 2] \\
10& [-4, -9, 30, 3] 	& 11 & [2, 4, 54, 3] 	&  12 & [-2, 13, 57, 3] \\
13 & [-3, -15, 42, 6] 	& 14 & [4, -9, 42, 6] 	&  15 & [0, -15, 45, 8] \\
16 & [-4, -8, 56, 8] 	& 17 & [3, -15, 48, 10] & 18 & [-3, -15, -48, -10]   \\
19 & [4, -8, -56, -8]  & 20 & [0, -15, -45, -8]  & 21 & [3, -15, -42, -6]  \\
22 & [-4, -9, -42, -6] & 23 & [2, 13, -57, -3]  &  24 & [-2, 4, -54, -3]  \\
25 & [4, -9, -30, -3]  &  26 & [-4, 16, 0, -2] & 27 & [4, 4, -28, -2]  \\
28 & [0, -3, -27, -2]  & 29 & [-4, -8, -16, -2] & 30 & [-2, 13, 3, -1] \\
31 & [2, 16, 14, -1]  & 32 & [2, -8, -14, -1]  & 33 & [-2, -11, -5, -1] \\
34 & [-4, -14, -2, -1] &     &                           &     & \\
\hline
\end{tabular}
\end{table}

The proof is a brute force approach. We check all points  $\p \in \wP_{(2,4,6,10)}$ of weighted height $\wh (\p) \leq 3/2$. 
There are, however, many rational points for weighted height   $2 < \wh (\p) \leq 3$.    ${\mathcal L}_2$ has at least these  rational points  as in \cref{tab:L2}.

\begin{lem}
There are no rational points $p\in {\mathcal L}_3$ with weighted  height $\wh (\p) < 2$. Moreover,   the only 
  the  rational points for weighted height   $2 < \wh (\p) \leq 3$ are displayed in \cref{tab:L3}
\end{lem}

\begin{table}[htp]
\caption{Rational points of height $\leq 3$ in ${\mathcal L}_3$}
\label{tab:L3}
\begin{tabular}{|c|c||c|c||c|c|}
\hline
 \# & $\p$   &   \# & $\p$ &   \# & $\p$  \\
\hline
1 &   [6, 18, 27, 2]  & 		2&   [-6, -18, 45, 2]  &  		3 &   [3, 18, 0, 4]  	\\
4& 	 [-3, -18, 36, 4] 	&		5 &	 [5, -26, 56, 4] & 		6 &	 [ -3, 27, 315, 4] \\
7 &	 [-5, 58, -76, 4 ] &		8 &	 [-5, 31, -49, 4 ] &  		9 &	 [ 5, 29, -9, 4 ]  \\
10 & [-2, -18, 39, 6]  &		11 & [-2, -18, 165, 6] &		12 & [2, 18, -15, 6] \\
13 & [-8, -80, 429, 8] &		14 & [-8, 49, -101, 8]  & 		15 & [8, -47, -59, 8]   \\
 16 & [-1, -18, 60, 12]   &		17 & [8, 36, 69, 12]   &      		18 & [-1, -45, 105, 12] \\
19 & [-8, -36, 123, 12]  &   	20 & [-5, 67, -55, 12] &		21 & [1, 18, -48, 12]  \\
 22 &  [1, 63, -15, 12] &		23 & [5, -65, -15, 12]  & 		24 & [6, 36, 36, 16]    \\
25&   [-6, -36, 108, 16] & 		26 &  [8, -63, 3, 24]    & 		27 & [-4, -36, 102, 24]  \\
 28 & [-8, 81, -171, 24]  &		29 & [4, 36, -6, 24]     & 		30 & [5, -59, 16, 32] \\
 31& [5, -68, 100, 32]   &   		32 & [-3, -36, 108, 32] &		33  & [-3, -27, 504, 32]  \\
34 &  [-5, 61, -132, 32] & 		35 & [3, 36, -36, 32]    & 		36  & [9, 54, 108, 36]  \\
37 & [-9, -54, 216, 36]& 		38 & [-2, -36, 132, 48] &            39  &  [-2, -72, 336, 48] \\
40 &   [3, 18, 0, 4]  &  		41 &     [-3, -18, 36, 4] &		42  &    [-2, -18, 39, 6] \\
43 &   [2, 18, -15, 6] &		44 &   [-1, -18, 60, 12] & 		45 & 					\\
\hline
\end{tabular}
\end{table}

\proof  The proof is a brute force approach. We check all points  $\p \in \wP_{(2,4,6,10)}$ of weighted height $\wh (\p) \leq 3$. 
\qed



\section{Determining if a genus two curve has split Jacobian using machine learning models}
Our next task will be to use machine learning to determine if a genus 2 curves has $(n, n)$ split Jacobian for relatively small $n$.    This has become an interesting problem to cryptographers lately; see \cite{costello-nn-split, nato-25, nato-elira}. 

Below we give an overview of what we obtained from running different machine learning models.  Throughout this section we assume the reader is familiar with basic concepts of artificial intelligence and machine learning in the level of \cite{AI}. 


\subsection{Preprocessing} 
Our dataset is structured as a dictionary containing 129 607 data points, with 30 399 belonging to $\mathcal{L}_3$ and 48 676 to $\mathcal{L}_2$. We begin by converting the dictionary into a pandas DataFrame (df2) to facilitate data manipulation. Tuple columns are then split into individual columns, and boolean columns are converted into binary values for compatibility with machine learning models.

Next, we select the columns 'J2', 'J4', 'J6', and 'J10' to construct our feature matrix X. Additionally, the columns 'In L3' and 'Automorphism group' are used as labels to evaluate the accuracy of the classification models. Before training any model, we normalize the data using the Normalizer scaler, which scales each sample to unit norm. This normalization ensures consistent scaling across samples, preventing features with larger magnitudes from disproportionately influencing the results.

\subsection{Unsupervised Learning}
First we will consider unsupervised learning.
\subsubsection{Autoencoder}
We assume the reader is familiar with autoencoders in the level of \cite{AI}. 
We implemented an autoencoder using TensorFlow and Keras to obtain a compact, latent representation of the data, which helps to capture important features and reduce dimensionality. The autoencoder consists of an encoder that maps the input data to a 3-dimensional latent space and a decoder that reconstructs the original data from this latent representation. The architecture includes dense layers with ReLU activation and Batch Normalization to stabilize training. To prevent overfitting, we used Dropout and early stopping with a patience of 5 epochs.

The autoencoder was trained on the normalized feature matrix, and the model was optimized using the Adam optimizer with a learning rate of 0.001 and mean squared error loss. By using this autoencoder, we aim to reduce the complexity of the data while preserving its most important features, making it easier for downstream models to analyze and interpret the data.
\subsubsection{K-means Clustering on latent representation.}
After obtaining the latent representation from the autoencoder, we applied KMeans clustering to the latent representation to identify distinct groups within the data. KMeans is an unsupervised algorithm that partitions the data into a specified number of clusters by minimizing the variance within each cluster. The best results was reached when the number of cluster was set to 4. 
\subsubsection{Gaussian Mixture Model (GMM) on latent representation}
In addition to K-Means, we applied the Gaussian Mixture Model (GMM) to the latent representation obtained from the autoencoder to model the data as a mixture of Gaussian distributions. GMM is a probabilistic model that assumes the data points are generated from a mixture of several Gaussian distributions, each with its own mean and variance. We set the number of components to 4 and used spherical covariance, which assumes that each Gaussian distribution has the same variance in all directions. After fitting the model, GMM predicted the cluster labels for each data point, allowing us to categorize the data. 

\subsubsection{K-means Clustering on original data}

We applied K-Means clustering directly on the original data to explore its clustering performance without using the latent representation obtained from the autoencoder.
The KMeans clustering on the original data gave a slightly less accurate results compared to clustering in the latent representation, demonstrating that the original features were less effective in distinguishing between the data groups.

\subsubsection{Gaussian Mixture Model (GMM) on original data}

We also investigated the performance of Gaussian Mixture Model (GMM) clustering on the
original data. GMM performed with almost the same accuracy  on the original data
compared to the latent representation obtained from the autoencoder. It achieved an
accuracy of 86\%, which was better compared to the KMeans model classification.

\subsubsection{Evaluation Measure}

To assess the quality of clustering, we employed the adjusted Rand score as our
evaluation metric. The adjusted Rand score measures the similarity between true
cluster assignments and the clustering results, accounting for chance agreement.
Higher scores indicate better clustering performance. We got these results:

\begin{itemize}
\item Adjusted Rand score for K-Means using autoencoder: 0.68
\item Adjusted Rand score for K-Means on original data: 0.65
\item Adjusted Rand score for GMM using autoencoder: 0.86
\item Adjusted Rand score for GMM on original data: 0.86
\end{itemize}

These \emph{Adjusted Random Index} (ARI) scores provide insights into the effectiveness of each clustering algorithm in
capturing meaningful patterns in the dataset. The higher ARI score for
GMM underscores its superior performance compared to K-Means.

\subsection{Supervised learning}
Supervised learning models exhibited superior performance in classifying curves into
the ${\mathcal L}_2$ and ${\mathcal L}_3$ spaces.
We explored several supervised learning models and found that both Random Forest and K-Nearest Neighbors (KNN) achieved the highest accuracy in predicting the clusters of genus 2 curves based on their invariant features. We also used a neural network which performed 96$\%$ accuracy in classifying ${\mathcal L}_2$ and ${\mathcal L}_3$ points.

For all the three supervised models we evaluated the accuracy using the test set and computed additional performance metrics like the F1 score to assess the balance between precision and recall. The F1 score is especially important when dealing with imbalanced datasets, as it gives a more balanced measure of the model's performance than accuracy alone.

We also plotted the learning curves to visually inspect the model's performance over time. These curves help identify issues like overfitting or underfitting by showing how the accuracy and loss change as the model trains on more data.

\subsubsection{Neural Network}
As a first step, the features X = ('J2', 'J4', 'J6', 'J10') were normalized using the Normalizer from Scikit-learn to ensure that all features have comparable scales. Then the target labels y were one-hot encoded to transform them into a format suitable for multi-class classification. 
The neural network consists of several layers:

   \textbf{Input Layer}: Input Layer has 4 units corresponding to the 4 invariant features.

   \textbf{Hidden Layers}: The model contains three fully connected (dense) layers with 128, 64, and 32 neurons, respectively. Each hidden layer uses the ReLU activation function to introduce non-linearity and improve the model's ability to capture complex patterns.

   \textbf{Batch Normalization}: After each dense layer, batch normalization is applied to stabilize learning and accelerate training by normalizing the activations.

  \textbf{Dropout}: Dropout layers with a rate of 0.4 were added after each hidden layer to prevent overfitting by randomly dropping a percentage of neurons during training.

   \textbf{Output Layer}: The output layer has as many neurons as there are classes, and it uses the softmax activation function to produce probability distributions over the classes.

During the model compilation and training step the model runs an Adam optimizer with a learning rate of 0.001 to adapt the learning rate during training. Also a categorical cross-entropy loss function was used for multi-class classification.
Early stopping was implemented with a patience of 10 epochs to stop training if the validation loss did not improve for 10 consecutive epochs. The model used ReduceLROnPlateau to reduce the learning rate by a factor of 0.2 if the validation loss stopped improving for 5 epochs.

The neural network model was evaluated using accuracy, classification report, and confusion matrix on the test data:

\emph{Accuracy}: The model achieved a test accuracy of 96$\%$, indicating the proportion of correctly classified points.

\emph{Classification Report}: The classification report provides additional insights into the precision, recall, and F1 score for each class, giving a better understanding of the model's performance on different classes; see \cref{tab:metrics}.

\begin{table}[h!]
\label{tab:metrics}
\caption{Classification Report for Neural Network model}
\centering
\begin{tabular}{lcccc}
\hline
Class & Precision & Recall & F1-Score & Support \\ \hline
1 & 0.87 & 1.00 & 0.93 & 6255 \\
2 & 1.00 & 0.90 & 0.95 & 9593 \\
3 & 1.00 & 1.00 & 1.00 & 10074 \\ \hline
Accuracy & \multicolumn{4}{c}{0.96 (25922 samples)} \\ \hline
Macro Avg & 0.96 & 0.97 & 0.96 & 25922 \\
Weighted Avg & 0.97 & 0.96 & 0.96 & 25922 \\ \hline   \\
\end{tabular}
\end{table}


\emph{Confusion Matrix}: The confusion matrix shows the model's classification performance across different classes, helping to identify any misclassifications.
\[
\begin{bmatrix}
6255 & 0   & 0    \\
938    & 8654 & 1  \\
2    & 0    & 10000
\end{bmatrix}
\]
The first row corresponds to the true class 1, with 6255 instances correctly classified as in ${\mathcal L}_3$, and no false positives or false negatives. The second row corresponds to the true class 2, with 8654 instances correctly classified as in ${\mathcal L}_2$, 938 instances of class 3 misclassified as class 2, and 1 instance of class 2 misclassified as class 1. The third row corresponds to the true class 3, with 10000 instances correctly classified as class 3, and 2 instances misclassified as class 1.

\emph{Learning curves}: The training and validation accuracy and loss were plotted over the course of training. These learning curves provide insights into the model's convergence, overfitting, or underfitting.
\begin{figure}[h!]
    \centering
    \includegraphics[width=0.75\textwidth]{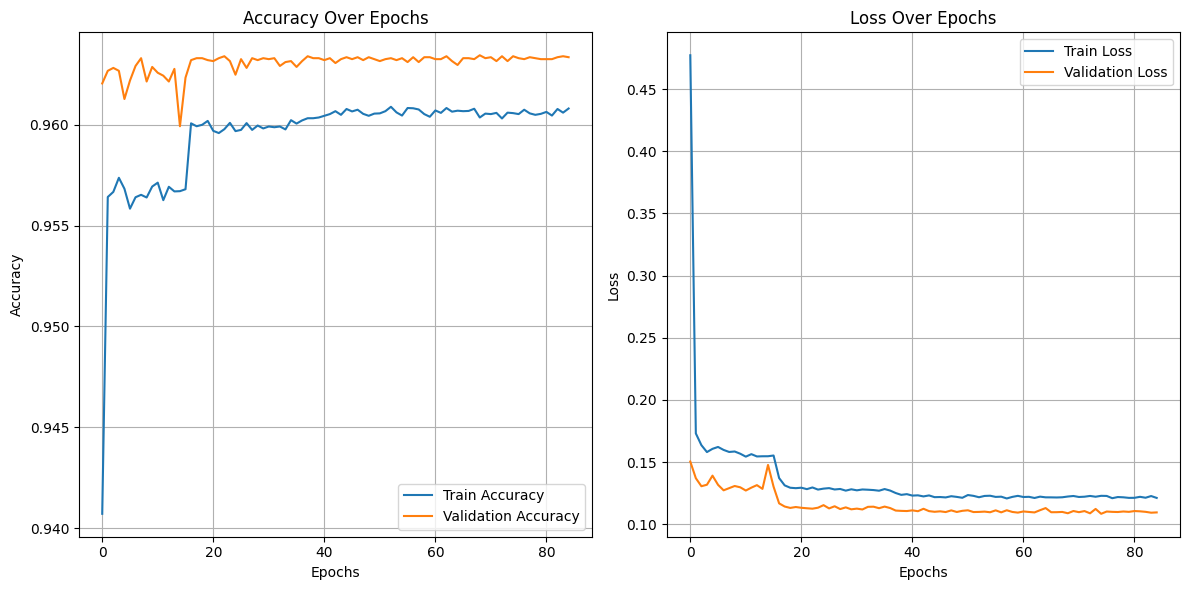}
    \caption{Learning curves for Neural Network model}
    \label{fig:example}
\end{figure}

\textbf{Accuracy}: The plot shows how the training and validation accuracy improved during training.

\textbf{Loss}: The loss plot shows the decrease in both training and validation loss, indicating that the model is learning effectively.
\subsubsection{Random Forest}
Random Forest is an ensemble learning method that constructs a collection of decision trees and combines their predictions. Each tree in the forest is built from a random subset of the data, and the final prediction is made by aggregating the predictions from all the trees (usually by majority vote for classification). This approach reduces overfitting and increases model robustness.

We trained the Random Forest model on the normalized training data with 200 estimators (trees) and evaluated it on the test set. The model achieved a notable accuracy score of 99.9$\%$, which indicates its effectiveness in distinguishing between the different clusters.
The confusion matrix for the classification results is as follows:
\[
\begin{bmatrix}
9315 & 0    & 0    \\
2    & 14513 & 2    \\
0    & 0    & 15051
\end{bmatrix}
\]
\begin{itemize}[leftmargin=*, itemsep=1mm]
\item The first row corresponds to the true class 1, with 9315 instances correctly classified as in ${\mathcal L}_3$, and no false positives or false negatives.
    \item The second row corresponds to the true class 2, with 14513 instances correctly classified as in ${\mathcal L}_2$, 2 instances of class 3 misclassified as class 2, and 2 instances of class 2 misclassified as class 1.
    \item The third row corresponds to the true class 3, with 15051 instances correctly classified as class 3, and no misclassifications.
\end{itemize}
The F1 score for the Random Forest model is 99.9\%. This indicates that the model performs exceptionally well, with both high precision and recall. The F1 score is the harmonic mean of precision and recall, which is calculated as follows:
\[
F1 = 2 \times \frac{\text{Precision} \times \text{Recall}}{\text{Precision} + \text{Recall}}
\]
A value of 99.9\% suggests that the model has very few false positives and false negatives, which is desirable in tasks where both types of errors have significant consequences.
\begin{figure}[hb]
    \centering
    \includegraphics[width=0.65\textwidth]{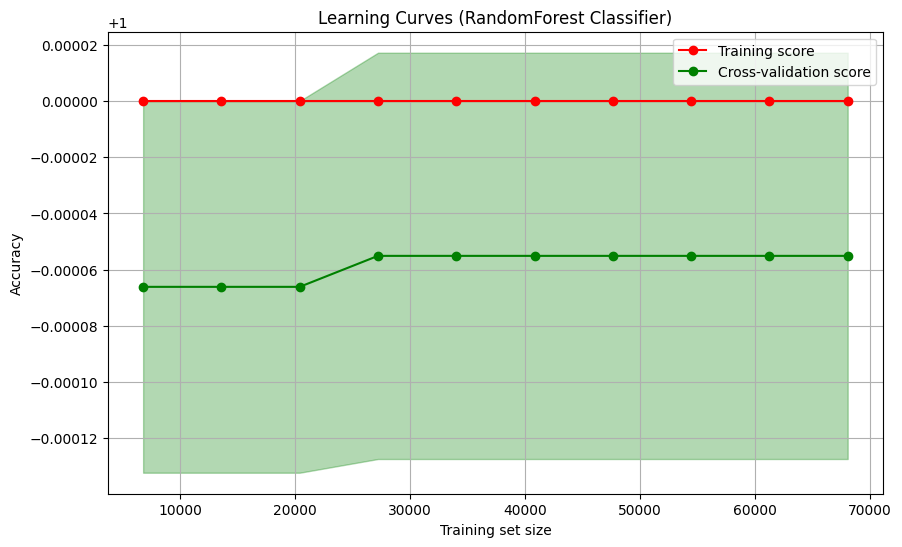}
    \caption{Learning curves for Random Forest model}
    \label{fig:rand-for}
\end{figure}

\subsubsection{K-Nearest Neighbors}
K-Nearest Neighbors (KNN) is a simple, instance-based learning algorithm that classifies a data point based on the majority class of its 'k' closest neighbors in the feature space. In this case, we used the Manhattan metric which calculates the distance between two points as the sum of the absolute differences of their coordinates. KNN is a non-parametric method, meaning it makes no assumptions about the underlying data distribution, which can be useful for datasets with unknown structures.

We applied KNN with 5 neighbors and trained it on the normalized data. The model achieved an accuracy score of 99.9$\%$, demonstrating its great performance on this task.
The confusion matrix for the K-Nearest Neighbors (KNN) model is as follows:
\[
\begin{bmatrix}
9312 & 3    & 0    \\
4    & 14511 & 2    \\
0    & 0    & 15051
\end{bmatrix}
\]
\begin{itemize}[leftmargin=*, itemsep=1mm]
    \item The first row corresponds to the true class 1, with 9312 instances correctly classified as in ${\mathcal L}_3$, and 3 instances in class 2, points in ${\mathcal L}_2$ misclassified as points in ${\mathcal L}_3$.
    \item The second row corresponds to the true class 2, with 14511 instances correctly classified as in ${\mathcal L}_3$, 4 instances of class 1 misclassified as class 2, and 2 instances of class 3 misclassified as class 2.
    \item The third row corresponds to the true class 3, with 15051 instances correctly classified as class 3, and no misclassifications.
\end{itemize}

Similar to the Random Forest model, the k-Nearest Neighbors (KNN) model achieves an exceptionally high F1 score 99.9\%, reflecting outstanding performance with very few misclassifications.

\begin{figure}[ht]
    \centering
    \includegraphics[width=0.8\textwidth]{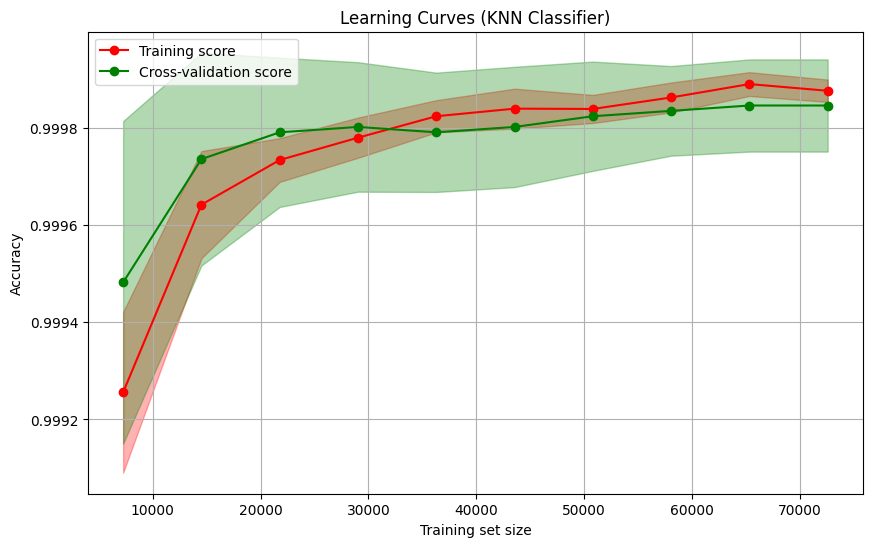}
    \caption{Learning curves for KNN model}
    \label{fig:example}
\end{figure}

The learning curves show a well-fitted model, The training and validation curves  both rise steadily and approach each other as training progresses. The gap between them remains small and steady.

 \section{Closing remarks}
The primary objective of this work was to explore the potential of artificial intelligence and data science techniques in addressing classical problems in arithmetic geometry, specifically focusing on the moduli space of genus two curves. While our main emphasis was on determining whether a moduli point has an $(n, n)$-split Jacobian—a topic of significant interest in isogeny-based cryptography—the methodologies developed here can be extended to a broad range of arithmetic questions related to genus two curves.

Current data science and machine learning methods predominantly operate on data within $\mathbb R^n$.  In contrast, we created a database of rational points in the weighted projective space $\mathbb P_{(2,4,6,10)}$.. An affine approach to this problem would require embedding  into a projective space, resulting in extremely large coordinates. Our choice of representation avoided this issue, and the results from our models demonstrated remarkably high accuracy. This high accuracy in detecting split Jacobians is notable, especially since we did not use the affine coordinates  of the data. Instead, the success of our models can be attributed to careful data normalization and the use of alternative distance metrics, such as Manhattan distance, rather than Euclidean distance.

Another significant factor contributing to the high accuracy is the sparsity of the data. As our results indicate, there are no curves with $(n,n)$-split Jacobians and small invariants (i.e., weighted height < 5). This natural sparsity, stemming from the way the data was generated, simplifies the model’s task of distinguishing between classes. The discovery of such distributions of split points is arguably the most important contribution of this work. Similar methods can be extended to higher values of $n$, such as $n=5, 7$, and beyond, provided one follows the techniques outlined in \cref{sec-4} and generates appropriate training data.

This study demonstrates the successful application of data science and machine learning techniques to problems defined over . Furthermore, by building upon our databases, these approaches can be extended to other arithmetic questions concerning genus two curves, such as those involving complex multiplication, endomorphism rings of Jacobians, isogenies, and more. These topics will form the basis of future research.

Our investigations also uncovered a new perspective: Can machine learning techniques be applied directly over graded vector spaces? In many scenarios, the input features of a machine learning model may be characterized by elements from a structured set . Neural networks defined over graded vector spaces could have significant mathematical and practical implications. For further discussion of such networks, we refer the reader to \cite{ann-graded-spaces}.

To the best of our knowledge, this work represents the first instance of using a machine learning model to study the geometry of weighted projective spaces. Remarkably, despite being trained on data consisting of rational points in a weighted projective space, the model appears to capture aspects of the Euclidean geometry of $\mathbb R^3$. This raises intriguing questions: Can a machine learning model be designed specifically for the geometry of weighted projective spaces ? Would such a model achieve higher accuracy and efficiency compared to models operating in $\mathbb R^n$? These questions will be addressed in forthcoming work \cite{ML-weighted, w-clusters}.


\section*{Declarations}
\subsection{Data availability} 
Data is available at  \href{https://www.risat.org/2024-03.html}{risat.org/2024-03}

\subsection{Funding}
There is no funding received from any funding agency. 

\subsection{Conflict of interest/Competing interests}
All authors declare that they have no conflict of interest to disclose.
 
\section*{Acknowledgements}
We want to thank the anonymous referees for very helpful observations and suggestion during the review process.

\bibliography{ref-2}

\end{document}